\documentclass[manyauthors]{Domfundam}
\pdfoutput=1
\usepackage{graphicx}
\usepackage{hyperref}
\usepackage{amsmath}
\usepackage{url}
\usepackage{amscd}
\usepackage{amsfonts}
\usepackage{amssymb}
\usepackage[]{babel}
\usepackage{pdflscape}
\usepackage{tikz}
\usepackage{tikz-cd}
\tikzcdset{scale cd/.style={every label/.append style={scale=#1},
		cells={nodes={scale=#1}}}}
\usetikzlibrary{decorations.pathmorphing}
\usepackage{todonotes}
\usepackage[T1]{fontenc}
\usepackage[utf8]{inputenc}
\usepackage{colortbl}
\usepackage{ulem}
\usepackage{dashbox}
\usepackage[top=2.2cm,left=2cm,right=1.4cm,bottom=1.5cm]{geometry}
\newtheorem{conjecture}{Conjecture}

\makeatletter
\newcommand{\vplus}{\mathbin{\mathpalette\vplus@\relax}}
\newcommand{\vplus@}[2]{%
	\raisebox{\depth}{\scalebox{1}[-1]{$\m@th#1\uplus$}}%
}
\makeatother

\usepackage{fancyhdr}

\begin{document}

	\pagestyle{fancy}
	\setcounter{page}{1} 
	\lhead{G. Mirkowska \& A. Salwicki}
	\chead{{\Large \thepage}}
	\rhead{Collatz conjecture becomes theorem /{\tiny \today}}
	\cfoot{}

	\title{Collatz conjecture becomes theorem \\{\small    $3n+1$ computations hide a paradox}}
	
	\author{Grażyna Mirkowska\\ UKSW University, Warsaw\\ Institute of Informatics  \and  Andrzej Salwicki \\Dombrova Research \\ Partyzantów 19,   05-092 Łomianki,    \textsc{Poland} \\ A.Salwicki \fbox{a} mimuw\textperiodcentered edu.pl }
	\date{\today}
	\bibliographystyle{alpha}
	\maketitle	
	\textbf{Keywords:} algorithm of Collatz,  halting property,  Collatz tree,  calculus of programs,algorithmic theory of numbers.  \bigskip\\
	ACM-class: F.3.1; D.2.4  \\
	MSC-class: 03D02 (Primary) 68Q02 (Secondary) \\
	
	\renewcommand{\thefootnote}{\arabic{footnote}} 	

	\begin{abstract} \color{white}I\color{black}\smallskip\\
					We are presenting the paradox, i.e. two theses T1 and T2 that contradict each other. Third thesis T3 solves the problem.\smallskip\\
			\ \ \textbf{T1}. 	\ \ 	We  show that the  Collatz conjecture  ''{For every natural number  $n$,  the $3n+1$    computation is finite.}'' is a \textit{semantically valid statement}.   \smallskip\\
			The sufficient  and necessary criterion $\varphi(n)$  for termination of  $3n+1$ computation  is given \ref{twefI}, c.f. page \pageref{cr3.1}.	\\	
			We prove that, every instance $\varphi (n/r)$ of the criterion where the variable $n$ is replaced by any natural number $r\neq0$, is a \textit{theorem} of Peano's arithmetic, Hence, the set $\left\lbrace  \varphi(n/r)\right\rbrace _{r=0}^{\infty} \subset Th(\mathcal{PA})$ is a recursive subset of the set of  theorems.  \smallskip\\
			\textbf{T2}. \ \    Paradoxically, the Collatz conjecture itself, \textit{is not a theorem} of number theory (Peano's arithmetic or a similar elementary theory).   \smallskip\\
			It is so because, 1\textdegree  the formula $\forall_{n}\varphi(n)$, obtained by putting the general quantifier $\forall_{n}$ in front of formula $\varphi(n)$,  may obtain  the value $\mathbf{\mathbb{F}}$ - false, in a  non-standard   model of Peano's arithmetic  \  and \ 
			2\,\textdegree there  is no way to bound the classical quantifier to the set of standard, reachable natural numbers.   \\ 
			To avoid the paradox, we will conduct our considerations in the formalized \textit{algorithmic} theory $\mathcal{ATN}$ of natural numbers.  The logical consequence operation of the theory is determined by  the   calculus of programs $\mathcal{AL}$, which is an extension of the predicate calculus.   \\
			The halting condition of the   computations  is expressed by  an algorithmic formula,   c.f. the consequent in  \eqref{Mnth1}.  \medskip\\
			\textbf{T3}. We are \underline{proving} that, four infinite sets $St_{0},St_{1},St_{2},St_{3}$ of formulas,  are the \textit{recursive sets} of \textit{theorems} of the theory $\mathcal{ATN}$.  Hence, every formula  of the set $St_{3}$  has a proof. Making use of the infinitary  inference rule $R_{3}$ (c.f.  page  \pageref{Rules}) to the infinite set $St_{3}$ of premises we conclude the proof of the Main theorem.
	\begin{scriptsize}
		\begin{equation}\label{Mnth1}\tag{{\tiny MainThm}}
			\mathcal{ATN} \vdash	\forall_{n>0}	\left( \underbrace{\left\{
			\begin{array}{l}
					q\leftarrow 1 ;  \\
					\mathbf{while}\ n \neq q \    \mathbf{do}\\ 
					\quad q\leftarrow q+1\\      \mathbf{od}
			\end{array}	
			\right\}(n=q)}_{{\mathbb{IF}\ n \ is\ a\ natural\ number}}	  \implies    
			\underbrace{\left\{ \begin{array}{l}
					m\leftarrow\rho(n);  \\
					\mathbf{while}\ m\neq 1 \    \mathbf{do}\\ 
					\quad m \leftarrow \rho(3m+1)\\      \mathbf{od}
				\end{array}
			\right\} (m=1)}_{\mathbb{THEN}\ the\ computation\   for\ n\ is\ finite\ \mathbb{FI}  } \right)  \qquad
		\end{equation}
	\end{scriptsize}
		\begin{scriptsize}  
	\hspace*{6.7cm}
		\textsc{Def.} The function $\rho$ is defined as  $\rho(n)=(2j+1) \stackrel{df}{\Leftrightarrow} \exists _{i}\exists_{j}\,n=2^{i}\cdot (2j+1)$.
	\end{scriptsize}
		\end{abstract}	
		
\maketitle
	\section{Introduction}\label{intro}

	The $3n+1$ problem remained open for over 80 years. It  has been formulated  in 1937 by Lothar Collatz, c.f. \cite{Lagar}.  The problem became quite popular due to its wording, for it is short and easy to comprehend.\\
	Collatz remarked that for any given natural number $n>0$, the sequence  $\{n_i\}$ defined by the following recurrence \ref{rec1}
	\[ \tag{rec1} \label{rec1} \left. \begin{array}{l}
		\nonumber m_0 = n  \\
		\nonumber m_{i+1} =  \begin{cases}  m_i  \div 2 &\mbox{  when $m_i$ is even}  \\
			3\cdot m_i+1 &\mbox{ when $m_i$ is odd}  \end{cases}  \quad\mbox{ for } \ i \geq 0
	\end{array}\right\}\] 
	seems always reach the value 1
	He formulated the following conjecture
	\[ \tag{Collatz conjecture} \text{for all }n\ \text{exists } i\ \text{such that } m_i=1 \]
	\textsc{Comment.}\ In the theory of recursive functions,  one would write down this problem as follows,  \\
	define the function $f_{C} \colon N \rightarrow N$ 
	\begin{equation*}
		f_{C}(n) =\begin{cases}
			f_{C}(n \div 2)  &  \text{when}\ n \ \text{is even} \\
			f_{C}(3n+1)  &  \text{when $n$ is odd and ,} \ n \neq 1  \\
			1 & \text{when }\ n=1
 		\end{cases}
	\end{equation*}
	and  ask, whether for every natural number $r$,  the value of the term  $f_{C}(r)$ is well defined.

	\subsection{ A bit of history}
Below, we  recall some important theorems, a couple of not too difficult remarks and a paradox.
	\begin{itemize}
		\item  The formulation of the Collatz problem contains the phrase "for every n in the set N of natural numbers".  Note, that the data structure of natural numbers can not be  axiomatized by any set of first-order formulas. It is a corollary of G{\"o}del's incompleteness theorem. Consquently,  no proof of Collatz conjecture exists in any first-order theory. 
		\item Note, that many important semantical   properties: "\textit{the Archimedean property}", "\textit{the halting property of program}", et al.,  can not be expressed  as first-order formulas \cite{ck:il, TMR}. 
		\item By the Kleene's theorem \cite{Kleen36}, the function defined bythe recurrence \ref{rec1} is calculated by  the following program \ref{Cl} \\
		\begin{equation}\label{Cl}\tag{Cl}
		 	\begin{tabular}{|l|} \hline
				\textbf{while} n $\neq$ 1 \textbf{do} \\
				\quad \textbf{if} even(n) \textbf{then}  n $\leftarrow$ n div 2 \
				 \textbf{else} n $\leftarrow$ 3n+1 \textbf{fi}  \\
				\textbf{od}   \\ \hline
			\end{tabular}
		\end{equation}
		\item   Erwin Engeler \cite{Eng:aps} in 1967 remarked that the halting property of program  can be expressed by an infinite disjunction of quantifier free formulas. \\
		Next, he proved that algorithmic properties of programs executed in the structure of ordered real numbers  are provable from the axioms of ordered,  fields plus the Archimedean property (understood as an infinite disjunction).
		\item  Programs are finite texts. The infinite disjunctions considered by Engeler have regularities. 
		\item Motivated by this observation   A. Salwicki (1969) introduced algorithmic formulas and the logical calculus (of programs) $\mathcal{AL}$ with iteration quantifiers \cite{al:gm:as}. The $\mathcal{AL}$ calculus allows to specify algorithmic properties  and to construct the proofs that verify the validity of properties. 
		\item  The following expression $\boxed{\left\lbrace \mathbf{while}\ \gamma\  \mathbf{do}\  \ x:=3x+1\ \mathbf{od}\  \right\rbrace \alpha(x)}  $ (it  consists of a program and a formula) is an example of simple algorithmic formula. Every pair consisting of a program and a formula is an algorithmic formula. The set of all formulas contains all classical formulas, all simple algorithmic formulas and is closed with respect to the boolean operations: disjunction, $\lor$ conjunction $\land$, negation $\lnot$ and two pairs of infinite boolean operations: the classical quantifiers $\exists. \forall$ and two iteration  quantifiers $\bigcap\left\lbrace K\right\rbrace \alpha$ to be read as \textit{the formula $\alpha$ holds after any iteration of the program $K$}. An existential iteration quantifier $\bigcup\left\lbrace K\right\rbrace \alpha$ reds:\textit{there is an  iteration of program $K$ such that the formula $\left\lbrace K^{i}\right\rbrace \alpha$ is satisfied.}  \\
	---------------------------  \\
		\item The logical value of the formula $K\alpha$ at a given valuation $v$ of variables  is determined as follows: \newline
		\begin{tabular}[b]{lp{6.7cm}}
				\begin{tikzcd} \\
				v \arrow[rr,"K\alpha"] \arrow[dr,"K"] &&\left\lbrace \mathbf{\mathbb{T},\mathbb{F} } \right\rbrace \\
				&v'\arrow[ur,"\alpha"] 
			\end{tikzcd}  & run the program $K$ at $v$ to arrive at another valuation $v'=K(v)$,  then evaluate the value of the formula $\alpha$ at the valuation $v'$ i.e.\qquad \ \  $val(K\alpha,v)=val(\alpha, v')$. \newline  If the result valuation $v'$ is not defined then $val(K\alpha,v)=\mathbf{\mathbb{F} } $.\\ \hline
		\end{tabular}
		\item   We introduced two countable families of infinite operations.   In addition to the classical quantifiers $\forall, \exists$, we consider the iteration quantifiers $\bigcup, \bigcap$.\\
		The meaning of the formula $\bigcap \left\lbrace K\right\rbrace\alpha $ is the greatest lower bound g.l.b. \ \  of the meanings  of the formulas of the set  $\left\lbrace \alpha,K\alpha, K^{2}\alpha,K^{3}\alpha,\dots\right\rbrace $.
		The meaning of the formula $\bigcup \left\lbrace K\right\rbrace\alpha $ is the least upper bound l.u.b.\  of the meanings  of the formulas  $\left\lbrace \alpha,K\alpha, K^{2}\alpha,K^{3}\alpha,\dots\right\rbrace  $. \\
		Algorithmic formulas of the form $\boxed{\mathbf{while}\ \gamma\  \mathbf{do}\ K\ \mathbf{od} \alpha}$  also define  infinite operations on logical values.	  
		\item The calculus of programs enjoys the \textit{completeness} property. \\
		 The proof of the completeness theorem   uses of the Rasiowa-Sikorski lemma, c.f. \cite{al:gm:as}. 
	\end{itemize}

		\subsection{An  illustration of some regularities.} 
	\quad  \\  Look at the  table \ref{tab1}. 
	
	\begin{table}[h]\begin{center}
				{\small\caption{	 	\textsc{Case} n=19  of (compact) computation} \label{tab1}
					\hspace*{1cm}		 1-st column  =$ m_{i} $  of   $Collatz $ computation,
					next columns:  x =no of multiplications, z=no divisions, \newline
					y=code of path from $ n $ to current $ m_{i} $.  
					In next columns 6 and 7 we check invariant $ n\cdot 3^{x}+y \stackrel{?}{=}m_{i}\cdot 2^{z}$.\smallskip\\	 }  	
			{\footnotesize 		\begin{tabular}{|c|c|} \hline
					\textbf{Computation} &  \textbf{Analysis}  \\ 
					$\begin{array}{r|rrr|rrrr}
						\textbf{m} & \textbf{x} & \textbf{y} & \textbf{z}& \textbf{k} & \mathbf{n3^{x}+y} &  m\cdot 2^{z}   \\\hline
						19 & 0 & 0 & 0 & 0 & 19& 19 \\
						29 & 1 & 1 & 1 &1 &  58 & 58  \\
						11   & 2 & 5 & 4 & 3 &176   & 176   \\
						17   &3 & 31 & 5 & 1 &544 & 544 \\
						13   & 4& 125& 7 & 2& 1664 & 1664 \\
						5    & 5& 503& 10&3& 5120 & 5120   \\
						1   &   6& 2533& 14& 4& 16384& 16384 \\ \hline
						1 & 7 &23983 & 16 &  2  & 2^{16}& 2^{16}  \\
						\cdots & &\cdots &\cdots & &\cdots &\cdots  
					\end{array}$  & 				
					\begin{tabular}{l}					
						For each $ i^{th} $-row \  
						$ 0 \leq i \leq 6 $ \\
						\begin{tabular}{l}
							$ k_{i+1} =\kappa(3m_{i}+1) $\\  $m_{i+1}=(3m_{i}+1)\div 2^{k_{i+1}}$\\
						$ x_{i}=i $ \ \quad \textsc{and} \ \  $ z_{i}=\sum_{l=0}^{i}k_{l} $  \\  					 
						\textsc{annd} $ y_{i}=\sum_{j=0}^{x_{i}-1} 3^{x_{i}-1-j}\cdot 2^{z_{j}} $ 
						\end{tabular}\\  
						 invariants:     
						 1\textdegree\  \   $n\cdot 3^{x_{i}}+y_{i} = m_{i}\cdot  2^{z_{i}}$ \\
						 2\textdegree\ $\ n\cdot \prod_{j=0}^{i-1}(3+\frac{1}{m_{j}})= m_{i}\cdot  2^{z_{i}} $  \\  
						two halting conditions: \ 1\textdegree\ $ \exists_i m_{i}=1 $ 
						\\  
				2\textdegree\  $\exists_{i \in N} \ n\cdot 3^{x_{i}}+y_{i} =   2^{z_{i}}$     \\						 
					\end{tabular}\\\hline
			\end{tabular}}			
	\end{center}
\end{table}
	An analysis of the algorithm \ref{Cl} leads to lengthy,  awkward  formulas that are error prone. \\
	However, it helped us in  noticing, that every computation of the 3n+1 algorithm is accompanied by a computation  on the triples $\left\langle  x,y,z\right\rangle $ . Here the value of $x$ shows how many operations of multiplications were done. Similarly, the value of $z$ is  the number of division by 2 executed. We remark that the value of $y$ is a code of path starting from $n$.  Look at table \ref{tab1}.\\
	The  computation on triples  always terminate. For some triples it ends with error. The problem reduces to the task of proving that for every number $ n $ there exists an approppriate, error-free triple.  \\
	
\subsection{Properties of the structure $\mathfrak{N}$ of the natural numbers.} \quad  \\
It is tacitly assumed that computations are performed in the standard structure of natural numbers	$\mathfrak{N}$.
\[\mathfrak{N} = \left\langle N, 0,1, s,+,*, =,odd\right\rangle   \]
where $N$ is the set $\{0,1, s(1), s(s(1)), \dots s^{p}(1), \dots\}$. \\
The meaning of the  functor $+$ is the standard operation of addition. 
The sign $=$ denotes the identity relation. The t predicate $odd(x)=\mathbf{\mathbb{T}}$ iff x is odd number.  Note, $odd(x)\stackrel{df}{\equiv} \exists_z x=z+z+1$\\
In fact, we do not need multiplication operation, for the term $3*x$  can be  replacded by $x+x+x$. Similarly $2*x=x+x$ and $x \div 2 =z$ iff $\exists_z x=z+z \lor x=z+z+1$. \\
\paragraph{A theory that completely specifies the structure $\mathfrak{N}$.}  \quad \\
We are going to formulate and prove a theorem of termination property of the Collatz algorithm. \\
Therefore we need to choose a theory that 1\textdegree) Specifies the structure $\mathcal{N}$. 2\textdegree) Allows to express the termination property of algorthm. 3\textdegree)  Provides enough tools  to  prove such formula.\\

Two traditional candidate theories are the  Peano's arithmetic or Presburger theory of addition (Remember our computa\-tions need not multiplication operation.).\\
However, we are using a third theory: the formalized,  algorithmic theory of natural numbers $\mathcal{ATN}$, \cite{al:gm:as,oct1:gm:as}. For it comes together with the following  \textit{categoricity property}. 

\begin{theorem}
	Every model of the $\mathcal{ATN} $ theory is isomorphic to  
	 the standard structure $\mathfrak{N}$. 
\end{theorem}

Note,  The operations of multiplication by 3  and division by 2  can be defined by  simple programs with addition as the only operation. This observation will be used in our considerations.

	\subsection{Properties of computations}\quad  \\
	Every natural number $ n $ may be presented in the form $ n=2^{i}\cdot(2j+1) $. C.f. the notion of the \textit{pairing function}.\\
	In the sequel we use two functions\qquad  $\kappa, \rho\colon N \rightarrow N$ \qquad defined as foolows
	\begin{definition}\label{df1.3}  \hspace*{0.5cm} 
		 \fbox{$\kappa(2^{i}\cdot(2j+1))  \stackrel{df}{=} i$} \hspace*{1.8cm}  \fbox{$\rho(2^{i}(2j+1)) \stackrel{df}{=}2j+1$}  
	\end{definition}
	 (In many texts the function $\kappa$ is written as $ exp(n,2) $.) \ Note that function $\kappa$ can be defined by an appropriate formula of Presburger arithmetic.  \\

\noindent		We shall present and analyse  computations in  their \textit{compact} form. The figure \ref{c19} illustrates the concept. 
		\begin{figure}[h]
			\begin{footnotesize}
				\begin{center}
							\begin{tikzcd} 
								&  {16}\arrow[ddd,near end, "/2^4"] &  40 \arrow[ddd,near end,"/2^3"] &  52 \arrow[ddd,near end,"/2^2"] &  34 \arrow[ddd,near end,"/2^1"] &88 \arrow[ddd,near end,"/2^3"]  & 58 \arrow[ddd,near end,"/2^1"]  & even\\
								\\
								\\
								& \rho(16) & \rho(40) \arrow[luuu,sloped,"3m+1"] & \rho(52)\arrow[luuu,sloped,"3m+1"] & \rho(34)\arrow[luuu,sloped,"3m+1"] & \rho(88)\arrow[luuu,sloped,"3m+1"] & \rho(58) \arrow[luuu,sloped,"3m+1"] & 19\arrow[luuu,sloped,"3m+1"]	& 	odd  \\
							compact		& 1    &5\arrow[l]& 13\arrow[l]  & 17\arrow[l] & 11\arrow[l] &29\arrow[l]& 19\arrow[l] 
							\end{tikzcd} 
							 \caption{An example of compact  computation} \label{c19}                                                                                    
		\end{center}   
		\end{footnotesize}		
					\end{figure}		
			\begin{definition}\label{defTmd}
				A \textit{compact} form of Collatz computation for a given number $ n $ is a subsequence that contains odd numbers only. More precisely, it is the sequence of numbers $ \left\lbrace  m_{i}\right\rbrace_{i \in I}  $ such that $ m_{0}={n}\div{2^{\kappa(n)}} $ and for $ i>0 $ the element $ m_{i+1}={(3m_{i}+1)}\div {2^{\kappa(3m_{i}+1)}} $, or simply $m_{i+1} = \rho(3m_{i}+1)$. 
			\end{definition}

\section{Collatz tree and halting criterion of $3n+1$ computations}\label{cotre}

			A universal (generic) definition of tree reads as follows
			\begin{definition}\label{dftree}
				A set $D$ of finite sequences of natural numbers is a \textit{tree} iff it satisfies the conditions (\textit{i}) and(\textit{ii}) 
				 \begin{itemize}				
					\item[(\textit{i})] if a sequence $s=i_{1},i_{2},\dots,i_{l-1},i_{l} \in D$ then its non-empty prefix $s'=i_{1},i_{2},\dots,i_{l-1} \in D$  also belongs to $D$,
						\item[(\textit{ii})] the empty  sequence $\emptyset \in D$. It is the root of the tree.
				\end{itemize}		
			\end{definition}		
			According to this definition the set of finite computations of the program \ref{Cl} written in the reverse order is a tree. We call it the Collatz tree.
The figure \ref{DC}	illustrates a fragment of the Collatz tree.

					\begin{conjecture}  \label{conj1}
						The Collatz tree contains all  natural numbers.
					\end{conjecture}
	
						\begin{figure}[h]
								\begin{center}						
						\resizebox{\textwidth}{!}{
							\begin{tikzpicture}
								[scale=0.61,transform shape,
								level 1/.style={sibling distance=30.2cm},
								level 2/.style={sibling distance=15.4cm},
								level 3/.style={sibling distance=13.5cm},
								level 4/.style={sibling distance=12.1cm},
								level 5/.style={sibling distance=11.7cm},
								level 6/.style={sibling distance=9cm},
								level 7/.style={sibling distance=7.3cm},
								level 8/.style={sibling distance=5.4cm},
								level 9/.style={sibling distance=3.7cm},
								level 10/.style={sibling distance=2.8cm},
								level 11/.style={sibling distance=2.0 cm},
								level 12/.style={sibling distance=1.1 cm},
								level distance=2cm
								]
								\tikzstyle{lien}=[->,>=stealth,rounded corners=5pt,thick ]
								\tikzset{individu/.style={ fill=#1!25},
									individu/.default={white}}
								\node [individu] {$  16$}[grow=up] 
								{ child { node [individu] {$  32$ } 
										child { node [individu]{64}  
											child { node [individu]{128}
												child { node [individu]{256}  
													child { node [individu]{512} 
														child { node [individu]{1024}  
															child { node [individu]{2048} 
																child { node [individu]{4096} 
																	child { node [individu]{8192} 
																		child { node [individu]{16384}  
																			child { node [individu]{32768} }
																			child { node [individu]{5461} }}
																	}
																	child { node [individu]{1365}
																		child { node [individu]{2730}
																			child { node [individu]{5460} }}
																	}
																}
															}
															child { node [individu]{341} 
																child { node [individu]{682}  
																	child { node [individu]{1364}
																		child { node [individu]{2728}  
																			child { node [individu]{5456} }
																			child { node [individu]{909}  }
																	} }
																	child { node [individu]{227}
																		child { node [individu]{454}  
																			child { node [individu]{908} }
																			child { node [individu]{151} }
																		}
																	}
																}
															}
														}
													}
													child { node [individu]{85}  
														child { node [individu]{170} 
															child { node [individu]{340}  
																child { node [individu]{680} 
																	child { node [individu]{1360}  
																		child { node [individu]{2720}
																			child { node [individu]{5440} }}
																		child { node [individu]{453}
																			child { node [individu]{906} }}               
																	}
																}
																child { node [individu]{113}
																	child { node [individu]{226}  
																		child { node [individu]{452}
																			child { node [individu]{904} } }
																		child { node [individu]{75} 
																			child { node [individu]{150} }}
																	}
																} 
															}
														}
													}
												}
											}
											child { node [individu]{$  21$}  
												child { node [individu]{42} 
													child { node [individu]{84} 
														child { node [individu]{168} 
															child { node [individu]{336} 
																child { node [individu]{672} 
																	child { node [individu]{1344} 
																		child { node [individu]{2688}
																			child { node [individu]{5376} } 
																	} }
																}
															}
														}
													}
												}
											}
									} }  
									child { node [individu] {$  5$} 
										child { node [individu]{$  10$ }  
											child { node [individu]{$  20$}
												child { node [individu]{40}  
													child { node [individu]{80} 
														child { node [individu]{160}  
															child { node [individu]{320} 
																child { node [individu]{640}  
																	child { node [individu]{1280} 
																		child { node [individu]{2560} 
																			child { node [individu]{5120} }
																			child { node [individu]{853} }
																		}
																	}
																	child { node [individu]{213} 
																		child { node [individu]{426} 
																			child { node [individu]{852} }
																	}}
																}
															}
															child { node [individu]{53}  
																child { node [individu]{106}  
																	child { node [individu]{212} 
																		child { node [individu]{424} 
																			child { node [individu]{848} }
																			child { node [individu]{141} }
																		}
																	}
																	child { node [individu]{35} 
																		child { node [individu]{70}  
																			child { node [individu]{140}} 
																			child { node [individu]{23}
																			}  
																	}}
																}
															}
														}
													}
													child { node [individu]{13} 
														child { node [individu]{26} 
															child { node [individu]{52} 
																child { node [individu]{104} 
																	child { node [individu]{208}  
																		child { node [individu]{416}
																			child { node [individu]{832}             
																		} }
																		child { node [individu]{69} 
																			child { node [individu]{138}             
																			}            
																		}
																	}
																}
																child { node [individu]{17} 
																	child { node [individu]{34}  
																		child { node [individu]{68}
																			child { node [individu]{136}             
																		} }
																		child { node [individu]{11} 
																			child { node [individu]{22}             
																			}            
																		}
																	}
																}
															}
														}
													}
												}
											}
											child { node [individu]{$  3$} 
												child { node [individu]{6} 
													child { node [individu]{12} 
														child { node [individu]{24} 
															child { node [individu]{48} 
																child { node [individu]{96}
																	child { node [individu]{192}
																		child { node [individu]{384}
																			child { node [individu]{768}
																			}
																		}
																	}
																}
															}
														}
													}
												}  
											} 
										}
									} 
								}
								;
						\end{tikzpicture}   }
						\caption{ A fragment of Collatz tree} 
					\begin{center}
							Shown levels 4-15.  It does not include levels 0-3,  they consist of elements  1  --- 2 --- 4  --- 8  --- . \\Does every natural number $n$ belong to the Collatz tree?
					\end{center}
						\label{DC}	
							\end{center}					
					\end{figure} 		
					  
					\subsection {Halting criterion for Collatz's computations.}]\label{4algos}			
					Our aim, in this section, is to find    a halting formula for the $3n+1$ computations that will exhibit the properties of the computations and therefore will facilitate  the proof. Note, that the obvious candidate,  i.e. the formula $\{Cl\}(m=1)$ leads to an infinite disjunction of quickly lengthening  (and obscure) conjunctions.  Therefore we shall consider other programs that are equivalent to the $Cl$ program.    \\ 
					We begin showing the   algorithm \eqref{Gr} equivalent to the algorithm $Cl$.\\
					 Next, we consider    three  algorithms $Gr1,
					Gr2,Gr3$ that are  successive extensions of the eq\ref{Gr} algorithm. \\
\paragraph*{Invariant of computation.}
					The algorithm $Gr2$ will allow to formulate an  \textit{invariant} of $3n+1$ computations.  \\
					The algorithm $Gr3$ will help to construct the necessary and sufficient condition for the termination of computations.
					\begin{lemma} \label{lemat5.1}The following algorithm \ref{Gr} 
						is equivalent to   Collatz algorithm $Cl$. 
					\end{lemma}   
							 \begin{scriptsize}
							 	\begin{equation}\label{Gr} \tag{Gr}
							 	\begin{array}{lll}
							 		\begin{array}{|l|ll}
									\hline 
									\textbf{while}\ even(n) \ \textbf{do} \ n:= n \div 2\ \textbf{od} ;\\ 
									\textbf{while}\ n \neq 1 \ \textbf{do}  \\
									\quad  n:= 3*n+1; \ \\
									\quad  \textbf{while}\ even(n) \ \textbf{do} \ n:= n \div 2\ \textbf{od}  \\
									\textbf{od} ; \\ \hline 
							\end{array} \qquad \mbox{\ \ or\ shorter}& &
								\begin{array}{|l|ll}
								\hline 
								n:=\rho(n); \\
								\textbf{while}\ n \neq 1 \ \textbf{do}  \\
								\quad  n:= 3*n+1; \ \\
								\quad n:=\rho(n);  \\
								\textbf{od} ; \\\hline 
							\end{array}   
						\end{array}
							 \end{equation}		
							 \end{scriptsize}	
					\begin{proof}
						The equivalence of the algorithms $Cl$ and $Gr$ is intuitive.  Compare the recurrence of Collatz (\ref{rec1}) and the following  recurrence (\ref{rec2} ) that is calculated by  the algorithm $Gr$.
						\begin{equation} \tag{rec2}\label{rec2}\left.
							\begin{array}{rcl}
								  k_0 = \kappa(n) &  \land \  & m_0 = \rho(n)   \\
								\nonumber   k_{i+1}  =  \kappa(3m_{i}+1)&\  \land   \  &  m_{i+1}  =  \rho(3m_i+1)     
								\qquad \text{for} \ i  \geq 0\end{array} \right\}   \end{equation} 
							The definitions of the functions $\kappa$ and $\rho$ are contained in Definition \ref{df1.3}, on page \pageref{df1.3}. \\
						Note, the following equivalence holds  
						 {\footnotesize  \begin{multline}
						  	\bigl(\kappa(n)=l' \land \rho(n)=m'\bigr)\Leftrightarrow  \ \left\lbrace\,l:=0;\,m:=n;\,\textbf{while}\ even(m) \ \textbf{do} \ m:= m \div 2;\,l:=l+1\ \textbf{od} \right\rbrace \bigl(m=m' \land l=l'\bigr)  . 
						  \end{multline}}
						One can say the algorithm $Gr$ is obtained by the elimination of \textbf{if} instruction from the $Cl$ algorithm. 
					\end{proof}  
					\begin{corollary}
						Every halting condition of the program $Gr$ is also a halting condition for the program $Cl$.
					\end{corollary}
					\paragraph*{Counting the operations of multiplication by 3 and division by 2.}
					Next,  we  present the algorithm $Gr1$, an extension of algorithm $Gr$. \smallskip\\
\begin{equation}\label{Gr1}\tag{Gr1}
						 \begin{array}{|l|}
									\hline 
									\textbf{var } n,aux,i :integer ;\ k,m :\textbf{arrayof }  integer; \\
									\quad \Gamma_1: \begin{array}{|l|} \hline  \ \ \ \ \  i:=0;\, k_i:=\kappa(n);
										\,  \,m_i:=\rho(n) ;\\  \hline
									\end{array}  \\
									\textbf{while}\ m_{i} \neq 1 \ \textbf{do}  \\
									\quad  \Delta_1: \begin{array} {|l|} \hline 
									  \quad   aux:= 3*m_{i}+1; \ \,
										  k_{i+1}:=\kappa(aux);\\ 
										\quad  m_{i+1}:=\rho(aux); \,
										  \,i:=i+1 \smallskip\\ \hline
									\end{array} \\
									\textbf{od}  \\\hline 
							\end{array}    
\end{equation}
					\begin{lemma}\label{lemma8}  Algorithm \ref{Gr1} has the following properties: \\
						(\textit{i})\ Algorithms \ref{Gr} and \ref{Gr1} are equivalent with respect to the halting property. \smallskip\\
						(\textit{ii})\ The sequences $\{ m_i\}$ and $\{k_i\}$ calculated by the algorithm \ref{Gr1} satisfy the recurrence \ref{rec2}.
					\end{lemma}
					\begin{proof}
						Both statements are very intuitive. Algorithm \ref{Gr1} is an extension of algorithm \ref{Gr}. The inserted instructions do not interfere with the halting property of algorithm \ref{Gr}.  Second part of the lemma follows easily from the remark that  $k_0=\kappa(\mathrm{n})$  and  $m_0=\rho(n)$ and  that for all $i>0$ we have  $k_{i+1}=\kappa(3*m_i+1)$  and $m_{i+1}=\rho(3\cdot m_i+1) $.  
					\end{proof}  
					Each odd number $m$, $m \in D$, initializes a new branch in the Collatz tree. Let us give a color number $x+1$ to each new  branch emanating from a branch with color number $x$. Note, for every natural number$p$ the set of branches of the color $p$ is infinite.
					Let $W_x$ denote the set of natural numbers that obtained the color $x$. The set $W_{x}$ will be called the $x$-th \textit{stratum} pf the set $N$ of natural numbers.\\

					\begin{footnotesize}
						\begin{figure}[h!]
							\resizebox{\textwidth}{!}{
								\begin{tikzpicture}
									[scale=0.5,transform shape,
									level 1/.style={sibling distance=16.2cm},
									level 2/.style={sibling distance=10.1cm},
									level 3/.style={sibling distance=7.5cm},
									level 4/.style={sibling distance=7.2cm},
									level 5/.style={sibling distance=7cm},
									level 6/.style={sibling distance=6cm},
									level 7/.style={sibling distance=3.6cm},
									level 8/.style={sibling distance=3cm},
									level 9/.style={sibling distance=2.1cm},
									level 10/.style={sibling distance=1.6cm},
									level 11/.style={sibling distance=1.3 cm},
									level distance=1.6cm
									]
									\tikzstyle{lien}=[->,>=stealth,rounded corners=5pt,thick ]
									\tikzset{individu/.style={ fill=#1!25},
										individu/.default={black}}
									\node [individu=green] {$  16$}[grow=up] 
									{ child { node [individu=green] {$  32$ } 
											child { node [individu=green]{64}  
												child { node [individu=green]{128}
													child { node [individu=green]{256}  
														child { node [individu=green]{512} 
															child { node [individu=green]{1024}  
																child { node [individu=green]{2048} 
																	child { node [individu=green]{4096} 
																		child { node [individu=green]{8192} 
																			child { node [individu=green]{16384}  
																				child { node [individu=green]{32768} }
																				child { node [individu=red]{5461} }}
																		}
																		child { node [individu=red]{1365}
																			child { node [individu=red]{2730}
																				child { node [individu=red]{5460} }}
																		}
																	}
																}
																child { node [individu=red]{341} 
																	child { node [individu=red]{682}  
																		child { node [individu=red]{1364}
																			child { node [individu=red]{2728}  
																				child { node [individu=red]{5456} }
																				child { node [individu ]{909}  }
																		} }
																		child { node [individu ]{227}
																			child { node [individu ]{454}  
																				child {node [individu]{908} }
																				child { node [individu=yellow]{151} }
																			}
																		}
																	}
																}
															}
														}
														child { node [individu=red]{85}  
															child { node [individu=red]{170} 
																child { node [individu=red]{340}  
																	child { node [individu=red]{680} 
																		child { node [individu=red]{1360}  
																			child { node [individu=red]{2720}
																				child { node [individu=red]{5440} }}
																			child { node [individu]{453}
																				child { node [individu]{906} }}               
																		}
																	}
																	child { node [individu]{113}
																		child { node [individu]{226}  
																			child { node [individu]{452}
																				child { node [individu]{904} } }
																			child { node [individu=yellow]{75} 
																				child { node [individu=yellow]{150} }}
																		}
																	} 
																}
															}
														}
													}
												}
												child { node [individu=red]{$  21$}  
													child { node [individu=red]{42} 
														child { node [individu=red]{84} 
															child { node [individu=red]{168} 
																child { node [individu=red]{336} 
																	child { node [individu=red]{672} 
																		child { node [individu=red]{1344} 
																			child { node [individu=red]{2688}
																				child { node [individu=red]{5376} } 
																		} }
																	}
																}
															}
														}
													}
												}
										} }  
										child { node [individu=red] {$  5$} 
											child { node [individu=red]{$  10$ }  
												child { node [individu=red]{$  20$}
													child { node [individu=red]{40}  
														child { node [individu=red]{80} 
															child { node [individu=red]{160}  
																child { node [individu=red]{320} 
																	child { node [individu=red]{640}  
																		child { node [individu=red]{1280} 
																			child { node [individu=red]{2560} 
																				child { node [individu=red]{5120} }
																				child { node [individu]{853} }
																			}
																		}
																		child { node [individu]{213} 
																			child { node [individu ]{426} 
																				child { node [individu ]{852} }
																		}}
																	}
																}
																child { node [individu]{53}  
																	child { node [individu]{106}  
																		child { node [individu]{212} 
																			child { node [individu]{424} 
																				child { node [individu]{848} }
																				child { node [individu=yellow]{141} }
																			}
																		}
																		child { node [individu=yellow]{35} 
																			child { node [individu=yellow]{70}  
																				child { node [individu=yellow]{140}} 
																				child { node [individu=blue]{23}
																				}  
																		}}
																	}
																}
															}
														}
														child { node [individu]{13} 
															child { node [individu]{26} 
																child { node [individu]{52} 
																	child { node [individu]{104} 
																		child { node [individu]{208}  
																			child { node [individu]{416}
																				child { node [individu]{832}             
																			} }
																			child { node [individu=yellow]{69} 
																				child { node [individu=yellow]{138}             
																				}            
																			}
																		}
																	}
																	child { node [individu=yellow]{17} 
																		child { node [individu=yellow]{34}  
																			child { node [individu=yellow]{68}
																				child { node [individu=yellow]{136}             
																			} }
																			child { node [individu=blue]{11} 
																				child { node [individu=blue]{22}             
																				}            
																			}
																		}
																	}
																}
															}
														}
													}
												}
												child { node [individu]{$  3$} 
													child { node [individu]{6} 
														child { node [individu]{12} 
															child { node [individu]{24} 
																child { node [individu]{48} 
																	child { node [individu]{96}
																		child { node [individu]{192}
																			child { node [individu]{384}
																				child { node [individu]{768}
																				}
																			}
																		}
																	}
																}
															}
														}
													}  
												} 
											}
										} 
									}
									;
							\end{tikzpicture}   }
							\caption{Strata  $W_{0} -- W_{4}$  of Collatz tree}   \label{fig2}
						\end{figure}
					\end{footnotesize}

					\noindent Let $s$ be a variable not occurring in algorithm $Gr1$.  The following lemma states the partial correctnes of the algorithm $Gr1$ w.r.t. precondition $s=n$ and post-condition $s \in W_i$. \\
					\begin{lemma} \label{lemat9}
						If the algorithm  $Gr1$ halts for a number $n$ then it computes the number $i$ of stratum $W_i$  that contains the number $n$, 
						\[\{Gr1\}(true) \implies \bigl((s=n) \implies \{Gr1\}(s\in W_i )\bigr)\]
					\end{lemma}
					\paragraph*{Invariant of computation.}
					Next, we present another algorithm \ref{Gr2} and a lemma \ref{propty}.  \smallskip\\ 
					\begin{equation}\label{Gr2}\tag{Gr2}
						\begin{array}{|l|}
							\hline 
							\textbf{var } n,aux,i,z,y :integer ;\ k,m :\textbf{arrayof }  integer; \\
							\quad \Gamma_2:\ \begin{array}{|l|} \hline   i,y:=0;\, z,k_i:=\kappa(n);
								\,  \,m_i:=\rho(n);\\  \hline
							\end{array}  \\
							\textbf{while}\ m_{i} \neq 1 \ \textbf{do}  \\
							\quad  \Delta_2: \begin{array} {|l|} \hline 
								\quad   aux:= 3*m_{i}+1; \ \,
								k_{i+1}:=\kappa(aux);\\ 
								\quad y:=3*y+2^{z_{i}};\,z:=z+k_{i+1};  \\
								\quad  m_{i+1}:=\rho(aux); \,
								\,i:=i+1 \smallskip\\ \hline
							\end{array} \\
							\textbf{od}  \\\hline 
						\end{array}    
					\end{equation}

					\begin{lemma} \label{propty}Algorithm $Gr2$ has the following properties: \\
						(\textit{i})  Both algorithms $Gr1$ and $Gr2$ are equivalent with respect to the halting property. \\
						(\textit{ii})  Formula $\varphi:\ \boxed{n\cdot3^i+y=m_i\cdot2^z}$ is  an \colorbox{green!!10}
						 {invariant} of the program $Gr2$ i.e. the formulas (\ref{Gamma2a}) and (\ref{Delta2a})
						\begin{equation}\label{Gamma2a} \{\Gamma_2\}\, \bigl( (n\cdot3^i+y=m_i\cdot2^z) \land i=0\bigr) \end{equation}    
						\begin{equation} \label{Delta2a} \bigl(( n\cdot3^i+y=m_i\cdot2^z) \land i=x\bigr)\implies \{\Delta_2\}\bigl(( n\cdot3^i+y=m_i\cdot2^z) \land i=x+1\bigr) \end{equation}
						are theorems of   the algorithmic theory of numbers    $\mathcal{ATN}  $.
					\end{lemma}
					\begin{proof}  One should conceive the formula $\varphi$ as an abbreviation of an lengthy and obscure expression. 
						Proofs of  formulas \eqref{Gamma2a}, \eqref{Delta2a} are easy, it suffices to apply the axiom of assignment instruction $Ax_{18}$, see subsection \ref{ATN}.
					\end{proof}
				
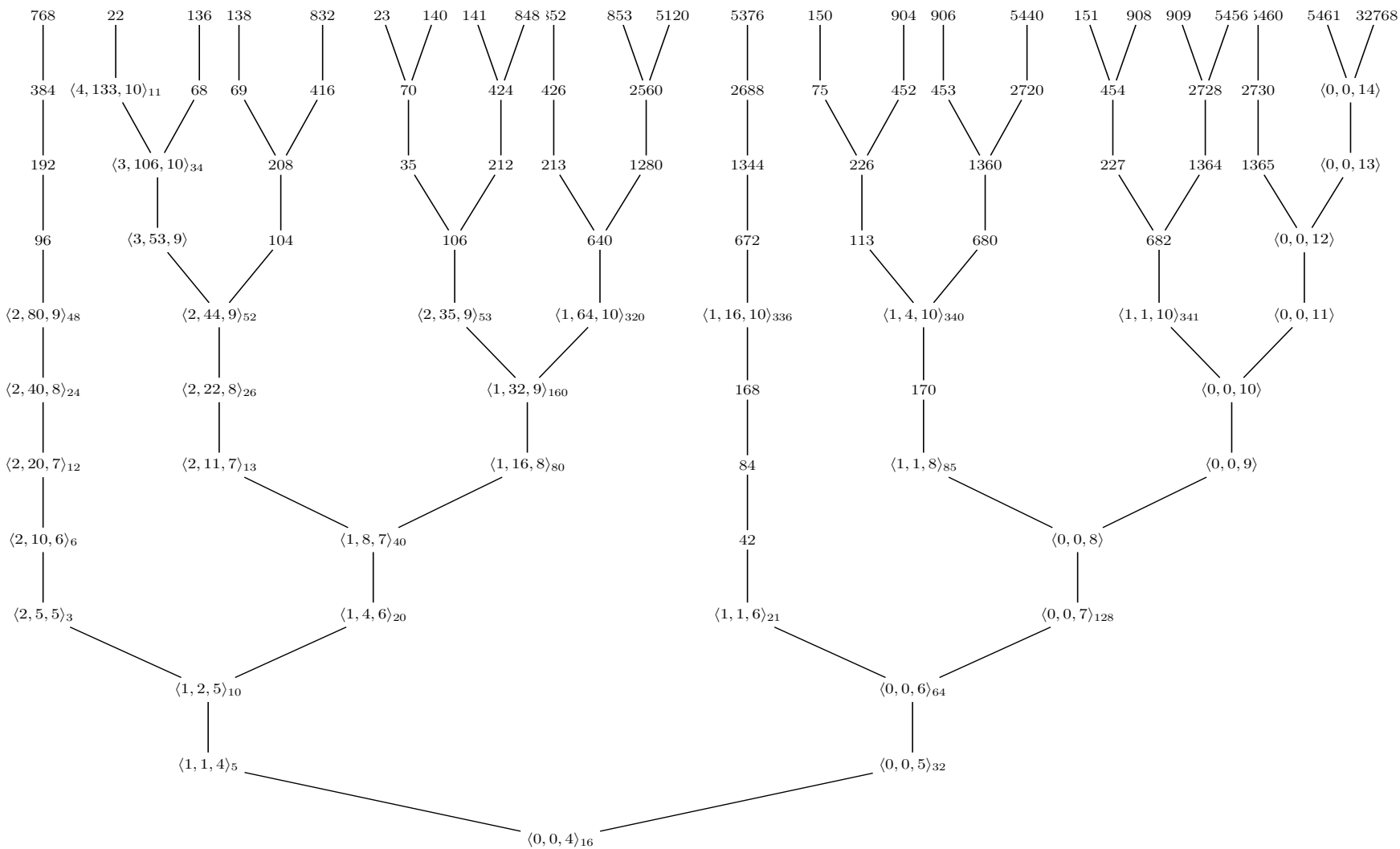
	\begin{figure}[h!]
		\resizebox{\textwidth}{!}{
	\begin{tikzpicture}
		[scale=0.48,transform shape,
		level 1/.style={sibling distance=16cm},
		level 2/.style={sibling distance=10cm},
		level 3/.style={sibling distance=7.5cm},
		level 4/.style={sibling distance=7.2cm},
		level 5/.style={sibling distance=7cm},
		level 6/.style={sibling distance=5.8cm},
		level 7/.style={sibling distance=3.3cm},
		level 8/.style={sibling distance=2.8cm},
		level 9/.style={sibling distance=2.1cm},
		level 10/.style={sibling distance=1.9cm},
		level 11/.style={sibling distance=1.2 cm},
		level distance=1.7cm
		]
		\tikzstyle{lien}=[->,>=stealth,rounded corners=5pt,thick]
		\tikzset{individu/.style={ fill=#1!25},
			individu/.default={black}}
		\node [individu=green] {$\langle 0,0,4 \rangle_{16}  $}[grow=up ] 
		{ child { node [individu=green] {$ \langle 0,0,5 \rangle_{32}  $ } 
				child { node [individu=green]{$\langle 0,0,6 \rangle_{64}  $}  
					child { node [individu=green]{$ \langle 0,0,7 \rangle_{128}   $}
						child { node [individu=green]{$\langle 0,0,8 \rangle  $}  
							child { node [individu=green]{$\langle 0,0,9 \rangle  $} 
								child { node [individu=green]{$\langle 0,0,10 \rangle  $}  
									child { node [individu=green]{$ \langle 0,0,11 \rangle  $} 
										child { node [individu=green]{$\langle 0,0,12 \rangle  $} 
											child { node [individu=green]{$\langle 0,0,13 \rangle  $} 
												child { node [individu=green]{$\langle 0,0,14 \rangle  $}  
													child { node [individu=green]{32768} }
													child { node [individu=red]{5461} }}
											}
											child { node [individu=red]{1365}
												child { node [individu=red]{2730}
													child { node [individu=red]{\ \ \ 5460} }}
											}
										}
									}
									child { node [individu=red]{$\langle 1,1,10 \rangle_{341}  $} 
										child { node [individu=red]{682}  
											child { node [individu=red]{1364}
												child { node [individu=red]{2728}  
													child { node [individu=red]{5456\ } }
													child { node [individu]{909}  }
											} }
											child { node [individu]{227}
												child { node [individu]{454}  
													child { node [individu]{908} }
													child { node [individu=yellow]{151} }
												}
											}
										}
									}
								}
							}
							child { node [individu=red]{$\langle 1,1,8 \rangle_{85}$}  
								child { node [individu=red]{170} 
									child { node [individu=red]{$\langle 1,4,10 \rangle_{340}  $}  
										child { node [individu=red]{680} 
											child { node [individu=red]{1360}  
												child { node [individu=red]{2720}
													child { node [individu]{5440} }}
												child { node [individu=red]{453}
													child { node [individu]{906} }}               
											}
										}
										child { node [individu]{113}
											child { node [individu]{226}  
												child { node [individu]{452}
													child { node [individu]{904} } }
												child { node [individu=yellow]{75} 
													child { node [individu=yellow]{150} }}
											}
										} 
									}
								}
							}
						}
					}
					child { node [individu=red]{$  \langle 1,1,6 \rangle_{21}$}  
						child { node [individu=red]{42} 
							child { node [individu=red]{84} 
								child { node [individu=red]{168} 
									child { node [individu=red]{$\langle 1,16,10 \rangle_{336}  $} 
										child { node [individu=red]{672} 
											child { node [individu=red]{1344} 
												child { node [individu=red]{2688}
													child { node [individu=red]{5376} } 
											} }
										}
									}
								}
							}
						}
					}
			} }  
			child { node [individu=red] {$\langle 1,1,4 \rangle_5  $} 
				child { node [individu=red]{$\langle 1,2,5 \rangle_{10}  $ }  
					child { node [individu=red]{$\langle 1,4,6 \rangle_{20}  $}
						child { node [individu=red]{$\langle 1,8,7 \rangle_{40}  $}  
							child { node [individu=red]{$\langle 1,16,8 \rangle_{80}  $} 
								child { node [individu=red]{$\langle 1,32,9 \rangle_{160}  $}  
									child { node [individu=red]{$\langle 1,64,10 \rangle_{320}  $} 
										child { node [individu=red]{640}  
											child { node [individu=red]{1280} 
												child { node [individu=red]{2560} 
													child { node [individu=red]{5120} }
													child { node [individu]{853} }
												}
											}
											child { node [individu]{213} 
												child { node [individu]{426} 
													child { node [individu]{852} }
											}}
										}
									}
									child { node [individu]{$\langle 2, 35,9 \rangle_{53}  $}  
										child { node [individu]{$\langle 2, 70,10 \rangle_{106}  $}  
											child { node [individu]{212} 
												child { node [individu]{424} 
													child { node [individu]{848} }
													child { node [individu=yellow]{141} }
												}
											}
											child { node [individu=yellow]{35} 
												child { node [individu=yellow]{70}  
													child { node [individu=yellow]{140}} 
													child { node [individu=blue]{23}
													}  
											}}
										}
									}
								}
							}
							child { node [individu]{$\langle 2,11,7 \rangle_{13} $} 
								child { node [individu]{$\langle 2,22,8 \rangle_{26} $} 
									child { node [individu]{$\langle 2,44,9 \rangle_{52} $} 
										child { node [individu]{104} 
											child { node [individu]{208}  
												child { node [individu]{416}
													child { node [individu]{832}             
												} }
												child { node [individu=yellow]{69} 
													child { node [individu=yellow]{138}             
													}            
												}
											}
										}
										child { node [individu=yellow]{$\langle 3,53,9 \rangle_{17}$} 
											child { node [individu=yellow]{$\langle 3,106,10 \rangle_{34} $}  
												child { node [individu=yellow]{68}
													child { node [individu=yellow]{136}             
												} }
												child { node [individu=blue]{$\langle 4,133,10 \rangle_{11} $} 
													child { node [individu=blue]{$\langle 4,266,11 \rangle_{22} $}            
													}            
												}
											}
										}
									}
								}
							}
						}
					}
					child { node [individu]{$\langle 2,5,5 \rangle_3 $} 
						child { node [individu]{$\langle 2,10,6 \rangle_6 $} 
							child { node [individu]{$\langle 2,20,7 \rangle_{12} $} 
								child { node [individu]{$\langle 2,40,8 \rangle_{24} $} 
									child { node [individu]{$\langle 2,80,9 \rangle_{48}$} 
										child { node [individu]{96}
											child { node [individu]{192}
												child { node [individu]{384}
													child { node [individu]{768}
													}
												}
											}
										}
									}
								}
							}
						}  
					} 
				}
			} 
		}
		;
\end{tikzpicture}    }
 	 
	\caption{	Tree of triples  (strata 4 -- 15)}
    \label{tDC}
	\end{figure}

\paragraph*{Criterion of termination.}
					\noindent Subsequent algorithm $Gr3$  exposes the history of  the calculations of $x,y,z$.		
					\begin{equation}\label{Gr3}\tag{Gr3}
						\begin{array}{|l|}
							\hline 
							\textbf{var } n,aux,i :integer ;\ k,m,x,y,z :\textbf{arrayof }  integer; \\
							\quad \Gamma_3:\  \begin{array}{|l|} \hline   i,y_{i}:=0;\, z_{i},k_i:=\kappa(n);
								\,  \,m_i:=\rho(n);\\  \hline
							\end{array} \\
							\textbf{while}\ n\cdot3^i+y_{i} \neq 2^{z_{i}} \ \textbf{do}  \\
							\quad  \Delta_3: \begin{array} {|l|} \hline 
								   aux:= 3*m_{i}+1; \ \,
								k_{i+1}:=\kappa(aux);\\ 
								 y_{i+1}:=3*y_{i}+2^{z_{i}};\,z_{i+1}:=z_{i}+k_{i+1};  \\
								  m_{i+1}:=\rho(aux); \,
								i:=i+1; x_{i}:=i \smallskip\\ \hline
							\end{array} \\
							\textbf{od}  \\\hline 
						\end{array}    
					\end{equation}
					The algorithm $Gr3$ has a couple of interesting properties.
					\begin{lemma}\label{lematGr3}
						Some semantical properties of the program Gr3 are:
						\begin{enumerate}
							\item[(\textit{i})] 	Both algorithms $Gr2$ and $Gr3$ are equivalent with respect to the halting property. 
							\item[(\textit{ii})] 	For every element $n$ after each  $i$-th iteration of algoritm $Gr3$, the following four formulas are satisfied 
								\begin{equation}
									\begin{array}{lll} \hline
										\varphi:\  n\cdot 3^{i}+ y_i = m_i\cdot 2^{z_i } &\qquad \ \qquad &
										\xi: x_i = i  \\
										\zeta:  z_i =\sum\limits_{j=0}^i\,k_j  & \qquad  \ \qquad& 
										\upsilon: y_i =\sum\limits_{j=0}^{i-1}\biggl(3^{i-1-j}\cdot 2^{z_j} \biggr)   \\  
										\hline
									\end{array} 
								\end{equation} 
							and the sequences $\{m_i\}$and $\{k_i\}$ satisfy the equalities  of the recurrence  (\ref{rec2}). 						
							\item[(\textit{iii})] 	In other words, the following formula is valid in the structure $ \mathfrak{N} $ of natural numbers
								\begin{equation}
								\mathfrak{N} \models \Gamma_3 \  \  \ \bigcap \, \{ \mathbf{if}\ m_i \neq 1\ \mathbf{then} \ \Delta_3 \ \mathbf{fi}\} \,(\varphi \land \xi \land \zeta \land \upsilon)
							\end{equation}
							\item[(\textit{iv})] for every element $n$  the sequence of triples $\langle  x_{i},y_i,z_i \rangle$  calculated by algorithm $Gr3$  is  increasing, monotone  as it can be seen from the following formula. 
							\begin{equation}
								(m_{i} \neq 1\land k_{i}=\kappa(3m_{i}+1) ) \implies 
								\left\lbrace \Delta_{3}\right\rbrace (x_{i+1}=i+1 \land y_{i+1}=3y_{i}+2^{z_{i}} \land z_{i+1}=z_{i}+k_{i}  )
							\end{equation} 
							\item[(\textit{v})]	Hence,  for every element $n$  the sequence of triples $\langle  X_{i}\,(=i),Y_i,Z_i \rangle$  calculated by algorithm $Gr3$  is  increasing, monotone  as it can be seen from the following formula \eqref{Increase}. 
							\begin{equation}\label{Increase}
								(m_{i} \neq 1\land k_{i}=\kappa(3m_{i}+1) ) \implies 
								\left\lbrace \Delta_{3}\right\rbrace (X_{i+1}=i+1 \land Y_{i+1}=3Y_{i}+2^{Z_{i}} \land Z_{i+1}=Z_{i}+k_{i}  )
							\end{equation}							
							\item[(\textit{vi})] after each $i$-th  iteration of subprogram $ \Delta_{3} $ the equivalence \eqref{cor3.2} holds
								\begin{equation}\label{cor3.2}
									\mathfrak{N} \models \boxed{\left\lbrace \Gamma_{3} \right\rbrace ;\left\lbrace \Delta_{3} \right\rbrace^{i}  (m_{i}\neq 1 \Leftrightarrow n\cdot 3^{i}+y_{i} > 2^{z_{i}})}
								\end{equation}
						\end{enumerate} 					
					\end{lemma} 
	
					\begin{remark} \label{rem17} We can say \textit{informally} that the algorithm $Gr3$ performs as follow \smallskip \\
						\hspace*{1cm} \qquad \begin{tabular}{l} \hline
							$i:=0$; \\
							\textbf{while} $n \notin W_i$ \textbf{do}  $i:=i+1$ \textbf{od} \\ \hline
					\end{tabular} \end{remark}
					Let us note an interesting information on route in graph
					\begin{corollary}\label{cor3.3} For every natural number $ n $, for every natural number $ p\in N $ the triple $ \left\langle p,y_{p},z_{p}\right\rangle  $ computed by the program $ \left\lbrace \Gamma_{3} ; \mathbf{for} \ i:=1\ \mathbf{to}\ p\ \mathbf{do}\ \Delta_{3}\ \mathbf{od} \right\rbrace  $ determines a path from the number $ n $ to number $ m_{p} $   
						\begin{equation}
							\mathfrak{N}  \models \boxed{\left\lbrace \Gamma_{3} ; \mathbf{for} \ i:=1\ \mathbf{to}\ p\ \mathbf{do}\ \Delta_{3}\ \mathbf{od} \right\rbrace   (  n\cdot 3^{p}+y_{p} =m_{p} \cdot   2^{z_{p}})}
						\end{equation}
					\end{corollary}
					Observe the following criterion 
				\begin{lemma} [\textit{The sufficient and necessary criterion of  termination of $3n+1$ computations}]\label{cr3.1}  
					Let $r$ be a natural number (i.e. a numeral).  The  following conditions are equivalent
					\begin{itemize}
						\item [(\textit{i})] the $3n+1$ computation for the number $r$  is finite,
						\item [(\textit{ii})] the following sentence \ref{twef} is valid in $\mathfrak{N}$
							\begin{equation}\label{twef}	
							\exists_{x}\,
							\left(  r\cdot 3^{x}+y=2^{z}  \land
								 \left(	\begin{array}{l}
								\left( y=	\sum\limits_{j=0}^{x-1}\,3^{x-1-j}\cdot 2^{\sum_{p=0}^{j}k_{p}}  \right) \land    
								\left(z= \sum\limits_{p=0}^{x}k_{p}\right) \land \\
								 \ \biggl(k_{0}=\kappa(r)\land m_{0}=\rho(r)  \biggr) \land  \\ 
								 \bigwedge\limits_{l=0}^{x-1} \biggl(k_{l+1}=\kappa(3m_{l}+1) \land m_{l+1}=\rho(3m_{l}+1)\biggr)   \end{array}\right) 
							\right) 	.	  		
						\end{equation}
							\item [(\textit{iii})] the following sentence \ref{twefI} is valid in $\mathfrak{N}$
						\begin{small}
							\begin{equation}\label{twefI}	
							 \{x:=0\} \bigcup\{x:=x+1\}\,
						\left( 	\left( r\cdot 3^{x}+y=2^{z}\right)  \land \left(
							\begin{array}{l}
								\left( y=	\sum\limits_{j=0}^{x-1}\,3^{x-1-j}\cdot 2^{\sum_{p=0}^{j}k_{p}}  \right) \land    
								\left(z= \sum\limits_{p=0}^{x}k_{p}\right) \land \\
								 \biggl(k_{0}=\kappa(r)\land m_{0}=\rho(r)  \biggr) \land  \\ 
								 \bigwedge\limits_{l=0}^{x-1} \biggl(k_{l+1}=\kappa(3m_{l}+1) \land m_{l+1}=\rho(3m_{l}+1) \biggr) 
							\end{array}\right) \right) 	.\quad\quad	  		
						\end{equation}
						\end{small}
					\end{itemize}			
				\end{lemma}	
	\begin{proof}
		The proof is left to the reader.
	\end{proof}
	Perhaps you have noticed the difference between formulas \ref{twef}  and \ref{twefI}. \\
	Here we are offering a couple of \textit{informal} explanations. \\
	The formula \ref{twefI} is the correct criterion for the finiteness of 3n+1 computations.  The logical value of it says, \textit{informally}, there exists $x$ such that $x=0 \lor x=1 \lor x=2 \lor ...$ and the conjunction contained in big parentheses is satisfied by $x$.   \\
	On the other hand, 	the classical quantifier $\forall$ does not provide such comfort. It does not exclude  the possibility that the value of $x$ is non-standard. 
	The formula \ref{twef} should only appear together with the reachability axiom (R) in the subsection \ref{ATN}. \\

								The question: \textit{is there a natural number $n$ such that its computation is infinite?}  arises in a natural way.
								The following lemma \ref{tr1} brings a partial answer to this questions. \medskip\\
								\begin{lemma}\label{tr1}   
									Let   $n \in N$ be an odd, natural number.
									The condition (\textit{i}) implies the condition (\textit{ii})
									\begin{itemize}
										\item [(\textit{i})] The number  $n$ satisfies te equality $ n+1= 2^{p+1}\cdot (2j+1)  $ where $p,j$ are natural numbers,
										\item[(\textit{ii})] the sequence $n=m_{1}<m_{2}< \dots<  m_{p}$ of $p$ consecutive numbers in the compact computation is monotone, increasing, and the next computed number is smaller   $m_{p}>m_{p+1}$.
									\end{itemize}				
								\end{lemma}
								\begin{proof} For every $ j=0,\dots,p $, 
									by the definition of compact Collatz computation, c.f. page\pageref{defTmd}, the following equality holds  \[ m_{j+1}=\rho(3\cdot m_{j}+1)  \]
									(\textsc{Base}) Let us look at $ m_{1}=\rho(3\cdot m_{0}+1) $	.
									Note,  $ 3\cdot m_{0}+1=3\cdot (2^{p+1}\cdot x -1)+1 =3\cdot (2^{p+1}\cdot x) -2=2\cdot (3\cdot 2^{p}\cdot x -1) $ and the value of the term $  3\cdot (2^{p}\cdot x -1) $ is an odd number, because $ expo(2\cdot 3\cdot (2^{p}\cdot x -1),2)=1 $.
									Hence $ m_{1}= 3\cdot (2^{p}\cdot x -1)$. \\
									(\textsc{Induction step}) In the same way  we show that for every $ j=2,\dots,p $ the subsequent number $ m_{j}  $ satisfies the equality  \eqref{mj}
									\begin{equation}\label{mj}
										m_{j} =\rho(3^{j}\cdot (2^{p-j+1}\cdot x -1)+1)    
									\end{equation} 
									and is odd. \smallskip\\
									It is obvious that the sequence $ n=m_{0},m_{1}, \dots,m_{p} $ is increasing.  
									It is evident that $ m_{p+1} < m_{p}$ because  the number  of the odd number $ m_{p} $ is  even, for the number $ 3^{p+1} \cdot x  $ is odd. 
								\end{proof} 
								\begin{corollary}\label{tr2}
									From the lemma \ref{tr1} we learn the following fact:
									\begin{itemize}
										\item For any natural number $ q \in Nat$ there exists a compact computation that contains the sequence of length $ q $ of consecutive, increasing odd numbers  
									\end{itemize}	
								\end{corollary} 
\begin{landscape}	
	\subsubsection{Lemma on partial correctnes:\hspace*{2cm} \ $\boxed{\left\lbrace Gr1\right\rbrace (m=1) \implies  \left\lbrace Gr3\right\rbrace (n\cdot 3^{x}+y=2^{z})}$	}								
	\begin{figure}[h]
		
		\begin{center}	
				\begin{footnotesize}
					\begin{tikzcd} 												
						r \arrow[r,red,"P"'] \arrow[d]    \arrow[rrrrr, bend left=15,"Gr3"] \arrow[rrrrr, bend left=7,red,"Gr1"]& m_1\arrow[r,red, "P"']\arrow[d,"k_{1}"]  &m_2 \arrow[d,"k_{2}"]\arrow[r,red, "P"']   &
						m_3 \arrow[d,"k_{3}"]\arrow[r,red, "P"']   &
						\cdots\  \     m_{i_{r}-1}\arrow[r,red,"P"']\arrow[d,"k_{i_{r}-1}"]  &m_{i_{r}}=1  \arrow[d,"k_{i_{r}}"]    			 \arrow[rrrrr,red,"S_{i_{1}}(...S_{i_{r-2}}(S_{ i_{r-1}}(S_{ i_{r}}(1)))...)",bend left=29]\arrow[rrr,red,"S_{i_{r-2}}(S_{ i_{r-1}}(S_{ i_{r}}(1)))",bend left=29] \arrow[rr,red,"S_{ i_{r-1}}(S_{ i_{r}}(1))",bend left=18] \arrow[r,red,"S_{ i_{r}}(1)"] & m_{r-1} \arrow[r,red,"S_{ i_{r-1}}"] &m_{r-2} \arrow[r,red,"S_{ i_{r-2}}"]&m_{r-3}  \arrow[r,red,"S_{ i_{r-3}}"] & \cdots m_1 \arrow[r,red,"S_{ i_{1}}"] &  r&\\
						t_{0}=\left\langle 0,0,0\right\rangle \arrow[r,"\bar s_{i_{1}}"] \arrow[rr,bend right=9,"\bar s_{i_{2}}(\bar s_{i_{1}}(t_{0})) "']\arrow[rrr,bend right=18,"\bar s_{i_{3}}(\bar S_{i_{2}}(\bar s_{i_{1}}(t_{0}))) "']
						\arrow[rrrrr, bend right=25,sloped,"\bar s_{i_{r}}(\dots(\bar s_{i_{3}}(\bar s_{i_{2}}(\bar s_{i_{1}}(t_{0}))))  \dots) "] 
						& t_{1}=\left\langle 1,1,2i_{1}\right\rangle   \arrow[r,"\bar s_{i_{2}} "]&   t_{2}  \arrow[r,"\bar s_{i_{3}}"]&
						t3  \arrow[r,"\bar s_{i_{4}} "] &
						\cdots\ \    t_{i_{r}-1}   \arrow[r,"\bar s_{i_{r}} "]&   t_{i_{r}}=\left\langle i_{r}, y_{i_{r}},z_{i_{r}}\right\rangle\arrow[r,blue,"\bar P"']\arrow[u]   \arrow[rrrrr, bend right=20,blue,"IC"']  &\bar t_{i_{r}-1} \arrow[r,blue,"\bar P"']   \arrow[u] &\bar t_{i_{r}-2}   \arrow[r,blue,"\bar P"']   \arrow[u] &\bar t_{i_{r}-3}    \arrow[r,blue,"\bar P"']   \arrow[u] &				\cdots\ \     \ \  \bar  t_{i_{1}}    \arrow[u]\arrow[r,blue,"\bar P"'] &  \left\langle 0,0,0\right\rangle \arrow[u]   & 											
					\end{tikzcd} 
				\end{footnotesize}   	
			\caption{\textsc{{\footnotesize Two algorithms to find a path from 1 to $r$ or from $r$ to 1}}} \label{fig10}
		\end{center}   
	\end{figure} 						
	The numbers $k_{j}$ and $i_{j}$ are related as follows $k_{j}=2i_{j} -\delta(m_{j-1})$.  
\qquad 	\color{red} For the definition of $\delta(n)$ see page \pageref{GraphG}. \color{black}
	\begin{normalsize}	
		\begin{lemma} [On the  partial correctness of Gr3 program] \label{lem4.8}	Let $r\neq 0$ be any natural number.  
			The following conditions are equivalent 
			\begin{itemize}
				\item [(\textit{i})] The algorithm  \eqref{Cl} halts.
				\item [(\textit{ii})] The algorithm \eqref{Gr3} halts and thepost-condition   
				$r\cdot 3^{x}+y=2^{z}$ holds $ \{Gr2\}(r\cdot 3^{x}+y=2^{z}) $
				\item[(\textit{iii})]   	{\scriptsize \begin{equation*}\label{shortpr}
						\{Gr2\}(n\cdot 3^{x}+y=2^{z})
						\implies 
						\underbrace{ \left\lbrace \begin{array}{l}
								\begin{array}{l}
									x':=x;\quad y':=y;\quad
									\quad z':=z; \quad m:=r;\\
									\mathbf{while}\ x'+y'+z' \neq 0 \ \mathbf{do} \\
									\   \begin{array}{l}
										k:=\kappa(y') ; \
										m:=P(m); \ 
										x':=x'-1; \
										y':=(y'\stackrel{\cdot}{\_}3^{x'}) \div 2^{k}; \
										z':=z'-k
									\end{array}  \\
									\mathbf{od} \\
								\end{array}
							\end{array}\right\rbrace}_{IC}  ( x'+y'+z'=0) 	\hspace*{3cm}				
				\end{equation*}}
		The triple $\langle x,y,z\rangle $ computed by the program Gr2 permits the algorithm IC to calculate the path from 1  to $r$.
			\end{itemize}
		\end{lemma}  	
		\begin{proof} 
		The proof 	is left to te reader.
		\end{proof}								
	\end{normalsize}

	\newpage				  	
	\section{Hotel Collatz $ \mathcal{HC} $}\label{hcgf}							
 					  						Look at the figure \ref{hoC}. Imagine that the hotel $  \mathcal{HC}  $ consists of infinitely many towers. 
					  	For each tower, the number of the room  on the ground floor of the tower is  odd. There is one room on each floor. Its number is twice the number of the room located on the floor below.
					  	Let $n=2^i\cdot(2j+1)$ . It means that the room number $n$ is located in $ j^{th}$-tower  on the floor number $i$ .   Each tower is equipped with an elevator (shown as a green line). Moreover, each tower is connected to another by a staircaise that connects numbers $k=2j+1$ and $3k+1$.  This is shown as a red arrow $\color{red}\overrightarrow{\color{black}\langle k,3k+1 \rangle} $. 
					  	\begin{definition}[Hotel Collatz]\label{defHC}
					  		The graph $\mathcal{HC}=\langle V, E \rangle$ is defined as follows\\
					  		\hspace*{1cm} $V= N$ \qquad {\footnotesize i.e. the set of vertices is the set of standard, reachable, natural numbers} \\
					  		\hspace*{1cm} $E=\{ \color{green}\overrightarrow{\color{black}\langle k,p \rangle} \color{black}: \exists_p\,k=p+p\} \cup \{\color{red}\overrightarrow{\color{black}\langle k,3k+1 \rangle} \color{black}: \exists_p k=p+p+1\} $    \qquad {\footnotesize are edges of the graph }
					  	\end{definition}  
					  	
					  						\begin{flushright}					  		
					  			\textbf{Note.} Don't forget, our drawing  \ref{hoC}  is only a small fragment of the infinite HC structure. \\
					  			The figure \ref{hoC} shows only a few  red arrows.  We drew only those red arrows that fit entirely on a page.   
					  	\end{flushright}	
					  	It's easy to see that the sloping (red) edges of the $\mathcal{HC}$ graph change direction to the left or right. An odd-numbered edge with an odd-numbered start goes to the right.
					  	An even-numbered edge with an even-mumbered origin  goes to the left.\medskip 			  	
						\begin{figure}[h!]
								\resizebox{24.3cm}{8cm}{  
										\begin{tikzpicture} [arrows=<-]   
												\tikzstyle{algo}=[rectangle,draw,fill=yellow!20]
												\tikzstyle{hc}=[rectangle]
												\edef\q{8}
												\edef\g{58}  
												\foreach \x in {0,...,\g}{
													\draw [green, very thick] [<<-] (\x*2,-67.7)--(\x*2,-39);
												}    
												
												\node[rectangle,minimum height=28cm,
												text = olive, text opacity=1, fill opacity=0.65,
												fill = green] (r) at (0.,-53.5) {\ \ \ \ \ \  };
												
												\node[rectangle,minimum height=28cm,
												text = olive, text opacity=1, fill opacity=0.5,
												fill = red] (r) at (4.,-53.5) {\ \ \ \ \ \ \  };    
												\node[rectangle,minimum height=28cm,
												text = olive, text opacity=1, fill opacity=0.5,
												fill = red] (r) at (20,-53.5) {\ \ \ \ \ \ \  };
												\node[rectangle,minimum height=28 cm,
												text = olive, text opacity=1, fill opacity=0.5,
												fill = red] (r) at (84,-53.5) {\ \ \ \ \ \ \  };
												\node[rectangle,minimum height=28 cm,
												text = olive, text opacity=1, fill opacity=0.5,
												fill = gray!95] (r) at (2.,-53.5) {\ \ \ \ \ \ \  };    
												\node[rectangle,minimum height=28 cm,
												text = olive, text opacity=1, fill opacity=0.5,
												fill = gray!95] (r) at (12.,-53.5) {\ \ \ \ \ \ \  };
												\node[rectangle,minimum height=28 cm,
												text = olive, text opacity=1, fill opacity=0.5,
												fill = gray!95] (r) at (52.,-53.5) {\ \ \ \ \ \ \  };
												\node[rectangle,minimum height=28 cm,
												text = olive, text opacity=1, fill opacity=0.5,
												fill = gray!95] (r) at (112.,-53.5) {\ \ \ \ \ \ \  };
												\node[rectangle,minimum height=28 cm,
												text = olive, text opacity=1, fill opacity=0.5,
												fill = yellow] (r) at (16.,-53.5) {\ \ \ \ \ \ \  };    
												\node[rectangle,minimum height=28 cm,
												text = olive, text opacity=1, fill opacity=0.5,
												fill = yellow] (r) at (34.,-53.5) {\ \ \ \ \ \ \  };    
												\node[rectangle,minimum height=28 cm,
												text = olive, text opacity=1, fill opacity=0.5,
												fill = yellow] (r) at (68.,-53.5) {\ \ \ \ \ \ \  };        
												\node[rectangle,minimum height=28 cm,
												text = olive, text opacity=1, fill opacity=0.5,
												fill = brown] (r) at (10.,-53.5) {\ \ \ \ \ \ \  };    
												\node[rectangle,minimum height=28 cm,
												text = olive, text opacity=1, fill opacity=0.5,
												fill = brown] (r) at (22.,-53.5) {\ \ \ \ \ \ \  };    
												\node[rectangle,minimum height=28 cm,
												text = olive, text opacity=1, fill opacity=0.5,
												fill = blue] (r) at (6.,-53.5) {\ \ \ \ \ \ \  };  
												\node[rectangle,minimum height=28 cm,
												text = olive, text opacity=1, fill opacity=0.5,
												fill = blue] (r) at (28.,-53.5) {\ \ \ \ \ \ \  };    
												\node[rectangle,minimum height=28 cm,
												text = olive, text opacity=1, fill opacity=0.5,
												fill = magenta] (r) at (8.,-53.5) {\ \ \ \ \ \ \  }; 
												\foreach \y in {1,...,\q} {
													\foreach \x in {0,...,\g} {   
		
														\pgfmathtruncatemacro{\result}{(2^(\q-\y))*(2*\x+1)}
														\edef\t{\result}
														\pgfmathtruncatemacro{\result}{\x*2}
														\edef\l{\result}
														\pgfmathtruncatemacro{\result}{\y*4}
														\ifnum \t<17000
														\node (\t) at (\l,-\result-36) {\t};
														\fi

													}
												}
												\draw [very thick,red][<-] (10)--(3);
												\draw [very thick,red][<-] (16)--(5);
												\draw [very thick,red][<-] (22)--(7);
												\draw [very thick,red][<-] (28)--(9);
												\draw [very thick,red][<-] (34)--(11);
												\draw [very thick,red][<-] (40)--(13);
												\draw [very thick,red][<-] (46)--(15);
												\draw [very thick,red][<-] (52)--(17);
												\draw [very thick,red][<-] (58)--(19);
												\draw [very thick,red][<-] (64)--(21);
												\draw [very thick,red][<-] (70)--(23);
												\draw [very thick,red][<-] (76)--(25);
												\draw [very thick,red][<-] (82)--(27);
												\draw [very thick,red][<-] (88)--(29);
												\draw [very thick,red][<-] (94)--(31);
												\draw [very thick,red][<-] (100)--(33);
												\draw [very thick,red][<-] (106)--(35);
												\draw [very thick,red][<-] (112)--(37);
												\draw [very thick,red][<-] (118)--(39); 
												\draw [very thick,red][<-] (124)--(41); 
												\draw [very thick,red][<-] (130)--(43); 
												\draw [very thick,red][<-] (136)--(45); 
												\draw [very thick,red][<-] (142)--(47);                 
												\draw [very thick,red][<-] (160)--(53); 
												\draw [very thick,red][<-] (352)--(117);                  
												\draw [very thick,red][<-] (226)--(75); 
												\draw [very thick,red][<-] (232)--(77); 
												\draw [very thick,red][<-] (244)--(81);                  
												\draw [very thick,red][<-] (340)--(113); 
												\draw [very thick,red][<-] (184)--(61); 
											\end{tikzpicture}
										}          
										\caption{Hotel Collatz}\label{hoC}   
									\end{figure} 
									\end{landscape}
								
\paragraph{Is the graph $\mathcal{HC}$ a tree?}
Instead of proving that the  $Collatz$ tree  contains every natural number  we shall study another question: is the graph $\mathcal{HC}$ a tree ? An eventual positive answer solves the $ 3n+1$ problem. \\
This problem in turn boils down to the question: is the graph $\mathcal{G}$  defined by the set of  all odd numbers, c.f. definition\ref{GraphG},  a  subgraph of the graph $\mathcal{HC}$, a tree?\\
So, we need to prove that the graph $\mathcal{G}$ is acyclic and connected.\\

\begin{corollary}
	It suffices to analyze the set of all odd, natural numbers. \\
	Imagine the graph $\mathcal{HC}$ in three dimensional space. Every tower stands over an odd number $o$. 
	Any two odd numbers $o$ and $p$ are connected by an edge iff there exists a natural number $k$ such that the follwing equality holds $p*2^{k}=3*o+1 $.
\end{corollary}		
We are going to prove the following conjecture.						
								\begin{conjecture}  \label{conj2}
									The hotel Collatz is an infinite, connected, acyclic graph, i.e. it is a tree.  Number 1 is the root of the tree.
								\end{conjecture}
				
\subsection{The graph $\mathcal{G}$ of odd, natural numers}	
						
	We shall introduce the notions of successors and of predecessor of an odd, natural number.
	The odd numbers divisible by 3 have no successors. 
	We shall use the following function $\delta$ to abbreviate the formulas.
	For very odd number  not divisible by 3, we put		
	 $\delta(n) \stackrel{df}{=} (n \mod 3) -1$.
	\begin{definition}\label{succ}	[of successor \quad  $m=S_{i}(n)$]  \newline
		The $i$-th \underline{successor} of an odd number $n$ is defined as follows 
		\[ S_{i}(n)\stackrel{df}{=} \dfrac{n\cdot 2^{2i-\delta(n)}-1} {3},\ \quad \    \text{where}\ \   i=1,2,\dots \ . \]
		\hspace*{7cm} \textsc{Note}, the odd numbers divisible by 3 ($n \mod 3 =0$) have no successors.  \\ 
		\hspace*{7cm} \textsc{Note} successor  $S_{1}(1)$  is not defined. 
	\end{definition} 
	\begin{definition}[of predecessor $n=P(m)$] The   \underline{predecessor}   of an odd number $m$	is the  odd number \ \  $n=\dfrac{3*m+1}{2^{\kappa(3m+1)}} $.  
	\end{definition}					
				    
We define the graph $\mathcal{G}$, To see an example subgraph of it, look at Fig. \ref{hofC31}, 	
	\begin{definition}\label{GraphG}
Graph $\mathcal{G}$ of odd numbers is the system of two sets
\[ \mathcal{G} \stackrel{df}{=} \left\lbrace V,E\right\rbrace \]	
The set $V$ (of \textit{nodes}) is the set of \underline{all} odd natural  numbers. \\
	The set $E$ (of \textit{edges }) is the set of all pairs $\left\langle n,m\right\rangle $ such that 
	\begin{enumerate}
			\item[c1)] $n,m$ are odd natural numbers, $E \subset V \times V$,
			\item[c2)] the number $n$ is indivisible by 3,\ \ $n\,\text{mod}\, 3 \neq 0$,
			\item[c3)] $m=\cdot \dfrac{n \cdot 2^{2i-\delta(n)}-1} {3}$  \quad  {\footnotesize where $i\in \{1,2,\dots\}$ and $\delta(n)\stackrel{df}{=}(n\,\text{mod}\,3) -1$}. 	
\end{enumerate}	
	\end{definition} 
	W may conceive the graph $\mathcal{G}$ as an algebraic system
	\[ \mathfrak{G} =\langle V; 1,\{S_{i}\}_{i=0}^{\infty},P,=\rangle \]
\noindent 
The following lemma gathers a couple of useful facts. 
\begin{lemma} \label{zbr}
	\begin{align}
		S_{i}(n)& \neq 1 \\
		S_{i}(n)& \neq n \\			
		  m=S_{i}(n)&\Rightarrow P(m)=n  \\
\label{rwn14}	 P(m)=n  &\Rightarrow  \exists_{i>0}\ m=S_{i}(n)&\\		 				  
		   S_{i}(n)=S_{j}(n)&\implies i=j \\
		   S_{i}(n)=S_{i}(m)& \implies n=m  \\
		     i \neq j &\implies S_{j}(S_{i}(n))\ \neq S_{i}(S_{j}(n)) \\
		     i <j &\implies S_{i}(n)\ < S_{j}(n) \\
		\text{the\ value\ of\ } S_{1}(1 )\ &\text{is\ undefined \ as well \ as\ the value\ of\  }P(1)\ &\  
		 	\end{align}
		\end{lemma}	
		\begin{proof}
			We shall prove the formula \eqref{rwn14}.  If $ P(m)=n$ then $ n\cdot 2^{\kappa(3\cdot m+1)}=3\cdot  m+1$. Hence $ i =(\kappa(3 \cdot m+1)+\delta(n))\div 2$. \\
		    The proof of the remaining formulas is left to the reader.
		\end{proof}

					\begin{figure}[h!]
			\begin{center}
				\begin{footnotesize}											
					\begin{tikzcd}  [row sep=0.7cm]
						&  	\vdots           \\
						\color{red}57	&&	& \\	
						43\dots\arrow[u,"2"]& & \color{red}177\dots  &\color{red}89\dots  \\									
						65\dots\arrow[u,"1"]&\color{red}33 & 133\dots\arrow[u]    & 67\dots\arrow[u] &\cdots \\	
						49\arrow[u,"2"]  & 197\dots &25\arrow[ul,"2"]\arrow[u,"4"] &101\dots\arrow[u]     & \color{red}51\dots  &205\dots&\cdots \\	
						\color{red}9& 37\dots\arrow[ul,"2 "]\arrow[u,"4 "]  &19\arrow[u,"2"]\arrow[ur,"4"] & 77\dots\arrow[ur,"1"]\arrow[urr,"3"]   & &325\dots   & \cdots \\	
						7\arrow[u,"2"]\arrow[ur,"4"] &29...\arrow[ur,"1"]\arrow[urr,"3"] &\color{red}15 & 61\dots\arrow[urr]    & \cdots \\		
						11\arrow[u,"1"]\arrow[ur,"3"] &\color{red}45&181\dots& \color{red}369 & 1477\dots& 23...\arrow[ulll,"1"]\arrow[ull,"3"']      &  \cdots \\	
						17\arrow[ur,"3"] \arrow[u,"1"]\arrow[urr,"5"] & \color{red}69  &  277...\arrow[ur,"2"] \arrow[urr,"4"]  &  35\arrow[urr,"1"]  &  \color{red}141...& 
						75  &301...  &\cdots \\
						\color{red}3 & 13\arrow[ul,"2"]\arrow[u,near end,"4"]\arrow[ur,"6"]& 53...\arrow[ur,"1"] \arrow[urr,"3"]&\color{red}213 &853\dots& 113\arrow[u,"2"] \arrow[ur,"4"]& \color{red}453\ \dots  & 227 &\color{red}909... 	 \\
						5\arrow[u,"1"',"S_{1} "]\arrow[ur,"3"'," S_{2}"]\arrow[urr,"5"]\arrow[urrr,"7"]\arrow[urrrr,near end,"9"',"S_{5} "]& &
						\color{red}21 &  85 \arrow[urr,"1"]\arrow[urrr,"3" ]  && 341\arrow[urr,near end,"1"]\arrow[urrr,near  end,"3"]   & \color{red}1365   && 5461...  \\	   
						&&1\arrow[ull,"4"',"S_{2}"] \arrow[u,"6"',"S_{3} "] \arrow[ur,"8"',"S_{4}"]\arrow[urrr,"S_{5}"] \arrow[urrrr,"S_{6}"] \arrow[urrrrrr,near end,"14"',"S_{7}"]    \\
					\end{tikzcd}   
					\caption{	{A fragment of the  graph $\mathcal{G}$, \scriptsize The\ numbers\ $k_{i}$\ on\ the \ edges  indicate \ the numbers ofsuccessors  $S_{i}$, namely  $i = \lceil \frac{k_{i}}{2} \rceil$}   } 					 
					\label{hofC31}							        
				\end{footnotesize}  
			\end{center}    
		\end{figure}  		
	One can easily remark, that the graph $\mathcal{G}$ is a forest, i.e. it is acyclic.
	 Is it a tree?			
\begin{lemma}\label{acyl}
	The graph  $\mathcal{G}$ is acyclic.
\end{lemma}
\begin{proof}
	The proof is easy, use the properties (13),(14),(17 - 21) listed in the lemma \ref{zbr}.
\end{proof}
It remains to be proved that the fraph $\mathcal{G}$ is connected.

		\subsection{The graph $\mathcal{G}$ is a  tree.}
	The following lemma  brings some crucial information.
\begin{lemma}\label{l3.3}
	For every element $n$ of the graph $\mathcal{G}$ such that $n \mod 3 \neq 0$  if there exist three natural numbers such that the following equality $n\cdot 3^{x}+y=2^{z} $ holds, then for every positive, natural number $i>0$ there exist three numbers $x_{i},y_{i},z_{i}$ such that the equality $S_{i}(n)\cdot 3^{x_{i}}+y_{i} =2^{z_{i}}$ holds.
\end{lemma}
{\footnotesize 	\begin{proof}
		\begin{align*}
			\intertext{From the assumption we have}
			n\cdot 3^{x}+y =& 2^{z} \\
			\intertext{Let $i>0$ be an odd natural numer. We multiply both sides of the equation  by $2^{2i-\delta(n)}$}
			(n\cdot 2^{2i-\delta(n)} )\cdot 3^{x}+(y\cdot 2^{2i-\delta(n)}) =& 2^{z} \cdot 2^{2i-\delta(n)} \\
			(n\cdot 2^{2i-\delta(n)} -1)\cdot 3^{x}+(y\cdot 2^{2i-\delta(n)}+3^{x}) =& 2^{z} \cdot 2^{2i-\delta(n)} \\
			\dfrac{n\cdot 2^{2i-\delta(n)} -1}{3}\cdot 3^{x+1}+(y\cdot 2^{2i-\delta(n)}+3^{x}) =& 2^{z} \cdot 2^{2i-\delta(n)} \\
			\intertext{We are putting 
				$x_{i}=x+1\ \& \ y_{i}= (y\cdot 2^{2i-\delta(n)}+3^{x})\ \& \ z_{i}=z+2i-\delta(n)   $  and obtain}
			S_{i}(n) \cdot 3^{x_{i}} + y_{i} = 2^{z_{i}}
		\end{align*} 
\end{proof}  }

We are going to address another interesting question: \textit{is it possible to permute the columns of the Hotel Collatz in such a way that all red edges go to the left?} (and therefore the the  odering $\prec$ of odd numbers  determined by this permutation of columns,  has the following property 
	$\forall_{n \in Odd}\  3 \nmid n \implies \forall_{i} \,  n \prec S_{i}(n)$ ?)
\begin{definition}(of subtrees $D_{c}$)
	For every ntural number $c$ we define a finite  tree $D_{c}=\langle V_{c},E_{c} \rangle$, a subtree of the graph $\mathcal{G}$.
	Let $p(o)$ be the length of the path from the odd number $o$ to 1.
	For every natural number $c$ we put $V_{c}=\{n \in Odd \colon p(n)\leq c \}$. and $E_{c}=E \cap V_{c} \times V_{c}$
\end{definition}
\begin{example}{A few  subtrees contained in the graph $\mathcal{G}$} \label{ex4.1}
\begin{center}
	\begin{small}
		$\begin{array}{r|cr}\hline
			c & D_{c}=\langle V_{c},E_{c} \rangle & \text{numbers on diagonal}\ c \\ \hline
			0 &\langle \{1\},\emptyset \rangle & 1 \\
			5 & \langle  \{1,5 \} , \{(1,5) \}\rangle & 5 \\
			7 &\langle \{1,5,21,3 \},\{(1,5),(1,21),(5,3) \}  \rangle  &21,3  \\
			9 &\langle \{1,5,21,,85,113 \},\{(1,5),(1,21),(5,3),(1,85),(5,13) \}  \rangle & 85,13\\
			11  &\langle \{1,5,21,3,85,13,341,53 \},\{(1,5),(1,21),(5,3),(1,85),(5,13),(1,341),(5,53) \}  \rangle & 341,53\\
			12 &\langle \{1,5,\dots,17,113\},\{(1,5),\dots,(85,113),(13,17)\}  \rangle & 13,17  \\			\hline
		\end{array}$ 
	\end{small} 
\end{center}
The trees form the increasing sequence $D_{0} \subset D_{5} \subset D_{9} \subset D_{11} \subset D_{12} \subset D_{13} \subset D_{14}$. \\
Every tree is bounded by a diagonal $c=x+z$. \\
The figure \ref{fig8}   illustraes the   trees $D_{0} -D_{14}$. \medskip\\
	\begin{figure}[h!]
\begin{tabular}[b]{p{\textwidth/2-0.8cm}|p{\textwidth/2-0.8cm}}\hline
	    {\tiny 		\begin{tikzcd}[,arrows=<-]
		& &&&&&&&& D_{0}\\
		1_{\langle 0,0,0 \rangle} 
\end{tikzcd}  }
&  	 {\tiny 		\begin{tikzcd}  [,arrows=<-]
		&  &&& &&&& D_{5}\\
		5_{\langle,1,1,4\rangle}   \\
		1_{\langle 0,0,0 \rangle}    \arrow[u]
\end{tikzcd}  } 

\\ \hline
	    {\tiny 		\begin{tikzcd}[,arrows=<-]
			    		3_{\langle 2,5,5 \rangle} & &&&&&&&& D_{7}\\
			    	5_{\langle,1,1,4\rangle} \arrow[u]&&21_{\langle 1,1,6 \rangle}\arrow[dotted,thick,green,ull] &&  &  \\
			    		1_{\langle 0,0,0 \rangle} \arrow[u]\arrow[urr] 
			    \end{tikzcd}  }
	   &  	 {\footnotesize 		\begin{tikzcd}[,arrows=<-]
	   		3_{\langle 2,5,5 \rangle}  & 13_{\langle 2,11,7 \rangle} && & D_{9}\\
	   		5 \arrow[u]\arrow[ur]&&21_{\langle 1,1,6 \rangle}\arrow[dotted,thick,green,ull] && 85_{\langle 1,1,8 \rangle} \arrow[brown,thick,dotted,brown,ulll]&  \\
	   		1 \arrow[u]\arrow[urr] \arrow[urrrr] 
	   \end{tikzcd}  } 

	   \\ \hline
	
	{\tiny 		\begin{tikzcd}[,arrows=<-,row sep = 1.0cm,column sep=0.35cm]
			\color{green}	3 & 13_{\langle 2,11,7 \rangle} &&53_{\langle,2,35,9\rangle}&&  D_{11}\\
			5 \arrow[u]\arrow[ur]\arrow[urrr] &&\color{green}21\arrow[dotted,thick,green,ull] && 85_{\langle ,1,1,8\rangle}\ \   \arrow[brown,very thick,dotted,brown,ulll]& \ \ 341_{\langle ,1,1,10\rangle}\arrow[red,thick,dotted,brown,ull] & \\
			1 \arrow[u]\arrow[urr] \arrow[urrrr] \arrow[urrrrr]
	\end{tikzcd}  } 
	&
	{\tiny 		\begin{tikzcd}	[,arrows=<-,row sep = 0.9cm,column sep=0.37cm]	
			& 17_{\langle,3,53,9\rangle} & & 	&& D_{12}  \\
			\color{green}	3  & 13\arrow[u] &&53_{\langle,2,35,9\rangle}& 113_{\langle,2,259,10\rangle}\arrow[dash,dotted,very thick,violet,ulll]&\\
			5 \arrow[u]\arrow[ur]\arrow[urrr] &&\color{green}21\arrow[dash,dotted,very thick,green,ull] && 85 \arrow[u] \arrow[dash,brown,very thick,dotted,brown,ulll]&  341_{\langle,1,1,10\rangle}\arrow[dash,red,very thick,dotted,brown,ull] &  \\ 
			1\arrow[u]\arrow[urr] \arrow[urrrr] \arrow[urrrrr]
	\end{tikzcd}  } 
	\\ \hline
		{\tiny 		\begin{tikzcd}[,arrows=<-,row sep = 1.2cm,column sep=0.24cm]		
			& 17_{\langle,3,53,9\rangle} & & 	&& D_{13}  \\
			\color{green}	3  & 13\arrow[u] &&53& 113_{,2,259,10}\arrow[dash,dotted,very thick,violet,ulll]&\color{green}213_{\langle,2,?,11\rangle}\\
			5 \arrow[u]\arrow[ur]\arrow[urrr]\arrow[urrrrr]  &&\color{green}21\arrow[dash,dotted,very thick,green,ull] && 85 \arrow[u] \arrow[dash,brown,very thick,dotted,brown,ulll]& 341 \arrow[dash,red,very thick,dotted,brown,ull] & \color{green}1365_{\langle 1,1,12\rangle} \arrow[dash,dotted,very thick,violet,ul]\\
			1 \arrow[u]\arrow[urr] \arrow[urrrr] \arrow[urrrrr]\arrow[urrrrrr]
	\end{tikzcd}  } 
	&
	{\tiny 		\begin{tikzcd}[,arrows=<-,row sep = 0.9cm,column sep=0.20cm]		
			& 11_{\langle,4,133,10\rangle} 
			&\text{Note!}& 35_{ 3,131,11} &\prec &69_{ 3,186,11} &D_{14}  \\
			& 17\arrow[u]&69\color{green}\arrow[dash,dotted,very thick,orange!90,ul]  &35\arrow[dash,dotted,very thick,orange!90,l]  & &&  \\			
		\color{green}	3  & 13\arrow[u] \arrow[ur]  &&53\arrow[u]& 113\arrow[dash,dotted,very thick,violet,ulll]&\color{green}213 &  227_{\langle,2,1027,12\rangle}\arrow[dash,dotted,very thick,orange!90,ulll] \\
	5 \arrow[u]\arrow[ur]\arrow[urrr]\arrow[urrrrr]  &&\color{green}21\arrow[dash,dotted,very thick,green,ull] && 85 \arrow[u] \arrow[dash,brown,very thick,dotted,brown,ulll]& 341 \arrow[ur]\arrow[dash,red,very thick,dotted,brown,ull] & \color{green}1365 \arrow[dash,dotted,very thick,violet,ul]\\
	1 \arrow[u]\arrow[urr] \arrow[urrrr] \arrow[urrrrr]\arrow[urrrrrr]  
\end{tikzcd} }

 \\ \hline
\end{tabular}  
\caption{Couple of trees $D_{c}$} \label{fig8}
	\end{figure}
Notice that each tree $D_{c}$,  is the lower-left part of the graph $G$. 
For every natural number $c$ the tree $D_{c}$ contains all odd numbers that are not farther than $c$ steps from 1. It means  that for every node $o \in V_{c}$ the computation has no more than $c=x+z$ steps and the equation $n\cdot 3^{x}+y=2^{z}$ holds.   
\end{example}
\begin{flushright}
End of 	Example  \ref{ex4.1} 	
\end{flushright}

\begin{landscape}	
	\begin{figure}
	\begin{tikzcd} [row sep = 0.9cm,arrows=<-]    
		&&&&  &&&&&&&&\color{white}X \arrow[dddd,red, dash,very thick] &\color{red}27\color{black}\arrow[dddd,blue,]  \\
		&&&&\color{red}57\arrow[d]  &&&&& &&\color{white}X \\
		&&&& 43\arrow[d]  &&&&& &\color{red}303\arrow[d,dotted]& 3077\arrow[dl,blue]\arrow[u,red,dotted] \\
		&&&&65\arrow[d] & \color{red}33\arrow[d]&&&& & 577\arrow[d,blue]  \\
		&&&& 49\arrow[d] &25\arrow[d] &&&&  & 433\arrow[d,blue] &&\color{white}\cdot\arrow[dddlll,red,dash,very thick] &\color{green}S_{1}^{17}\arrow[d,blue] \\
		& && \color{red}9\arrow[d] & 37\arrow[dl]  & 19\arrow[d] &&&&\color{red}81\arrow[d]  & 325\arrow[dl,blue] & & \color{red}111\arrow[d,dashed]&445\arrow[dl,blue]\\
		& && 7\arrow[d]  && 29\arrow[dll] & &&\color{red}15\arrow[d,dashed ] & 61\arrow[dl,blue]   && &\!\!\!\!\!\!\!\!\!167\color{green}{\footnotesize =S_{1}^{10}(2429)}\arrow[d,color=blue, thick]\\
		 & && 11\arrow[d] &&& \color{red}45\arrow[dlll] && 23\arrow[d,blue] &\color{white}X\arrow[dddd,red,dash,very thick]& &607\arrow[d,dashed]&2429\arrow[dl,blue]\\
		& && 17\arrow[d]  &&&& \color{red}69\arrow[dllll] & 35\arrow[d,blue] & && 911\arrow[d,blue]\\
		& \color{red}3\arrow[d] && 13\arrow[dll] &&&&& 53\arrow[dlllllll,blue]  & &&1367\arrow[d,blue] \\
		& 5\arrow[dl,blue]  & \color{red}21\arrow[dll] &&&&&& &\color{white}X& & 2051\arrow[d]\arrow[d,blue] \\
		1 &&&&&&&& &\color{white}X&& 3077 \\
	\end{tikzcd}     
	\caption{A subgraph of the graph $\mathcal{G}$ from another perspective\qquad 
		$ 27=S_{1}^{17}(S_{2}(S_{1}^{10}(S_{2}(S_{1}^{3}(3077)))))  \ \ \text{and}
		 \ \         3077=S_{2}(S_{1}^{2}(S_{2}^{2}(S_{1}^{2}(S_{3}({2}(1))))))$}
\end{figure}       
\end{landscape}

\begin{lemma}\label{l4.3}
	The trees $D_{c}$ have the following properties
	\begin{enumerate}
			\item [(\textit{i})] The trees $D_{0},D_{5},D_{7},D_{9},D_{11}$ are defined as shown in the example \ref{ex4.1}.
			\item [(\textit{ii})] Trees $D_{1},D_{2},D_{3},D_{4},D_{6},D_{8},D_{10}$ are not defined.
		\item [(\textit{iii})] For every $c \geq 11$ the following inclusion holds $D_{c} \varsubsetneq D_{c+1}$
		\begin{small}
			For $c+1=5,7,9,11$ trees $D_{c+1}$ are constructed on the base of trees $D_{1},D_{5},,D_{7},D_{9} $ respectively.
		\end{small}
		\item [(\textit{iv})] For every tree $D_{c}$ and  for every number $n \in D_{c}$ 
		if a successor $S_{i}(n)$ has a path of length $c+1$ then the number $n$ is in the tree $D_{c}$, i.e.  $n \in D_{c}$
		\item [(\textit{v})] For every odd number $n$ it exists the number $c$ such that the tree $D_{c}$ contains the initial segment $1,3, \dots, n$  of odd natural numbers.
	\end{enumerate} 
\end{lemma}

We shall order odd, natural numbers in a way similar to the rdering pairs of natural numbers defined by G. Cantor.  \\
\begin{definition}\label{precr}
	Let $a$ and $b$ be two odd, natural numbers. We  assume that the number $c$ is the least natural number that $a$ is in the tree $D_{c}$. Similarly, we assume that $b \in D_{d}$ and for every number $d'<d$ the relation $b \notin D_{d'}$ holds.    By the lemma \ref{l3.3} we know that exists two tiples $\langle x,y,z \rangle$ and $\langle u,v,t \rangle$ sch that two equations hold $a\cdot 3^{x}+y=2^{z}$ and $b\cdot 3^{u}+v=3^{t}$.    \\
	We put 
\begin{equation}
	a \prec b \stackrel{df}{\equiv} x+z <u+t \lor x+z=u+t \land x <u \lor x+z=u+t \land x =u \land y<v
\end{equation}
\begin{example}\begin{tabular}{|l|}
	 An initial segment of the sequence of natural numbers as ordered by the $\prec$ relation  \\
	$1\prec5\prec\color{green}21\prec3\color{brown}\prec85\prec13\prec\color{red}341\prec53\prec\color{blue}113\prec17\prec\color{violet}1365\prec213\prec
	\color{orange}227\prec35\prec69\prec11\prec\dots$
\end{tabular}
\end{example}
\end{definition}
To every odd natural number $n$ we assign its Collatz number  $C(n) $
\[ C\colon N \rightarrow N \]
\begin{definition}\label{Cnum}
C(n)$\stackrel{df}{=} \mathbf{card}(\{m \colon  	m \prec n  \})+1$
\end{definition}
\begin{landscape}
\begin{figure}

\begin{tikzcd}	 [column sep=0.55cm,row sep=0.6cm, arrows=<-]
	8192  \\
	4096 \arrow [u,green,thick]  \arrow [ddddddddddddrrrrrrrrrr] \\
	2048\arrow [u,green,thick]   \\
	1024\arrow [u,green,thick]  \arrow [ddddddddddrrrrrr]  \\
	512\arrow [u,green,thick]   \\
	256\arrow [u,green,thick] \arrow [ddddddddrrrr]  \\
	128 \arrow [u,green,thick] & 640 \arrow[dddddddrrrrrrrrrr]\\
	64\arrow [u,green,thick]  \arrow [rrdddddd]  & 320\arrow [u,red,thick]  \\
	32\arrow [u,green,thick]  &160\arrow [u,red,thick] \arrow [dddddrrrrrr]\\
	16\arrow [u,green,thick]  \arrow [rdddd]& 80\arrow [u,red,thick] &&&&208 \arrow [ddddrrrrrrrrr] \\
	8\arrow [u,green,thick]  &40\arrow [u,red,thick] \arrow [dddrrrr] &&&&104\arrow [u,thick] \\
	4\arrow [u,green,thick]   &20 \arrow [u,red,thick] \arrow [u]   && & 340\arrow [ddrrrr] & 52\arrow [u,thick]\arrow [ddrrrr]  \\
	2\arrow [u,green,thick]   & 10\arrow [u,red,thick] \arrow [rrd,dashed]&   && 170\arrow [u,thick] &26\arrow [u,thick] &682\arrow [drrrrrr,dashed]  & 106\arrow [drrrrrr,dashed] && 34\arrow [drrrrrr,dashed] &&&& \\
		1\arrow [u,green,thick] &5\arrow [u,red,thick] \arrow [u]  &\color{green}21&\color{green}3&\color{brown}85\arrow [u,thick]&\color{brown}13\arrow [u,thick] &\color{red}341\arrow [u,thick]&\color{red}53\arrow [u,thick]&\color{blue}113&\color{blue}17\arrow [u,thick]&\color{violet}1365&\color{violet}213&
	\color{orange}227&\color{orange}35&\color{orange}69&\color{orange}11&\dots  \\
\end{tikzcd}
\caption{Hotel Collatz with columns ordered by the relation $\prec$. All diagonal arrows point to the left.} \label{HC2}
\end{figure}	

\end{landscape}
\begin{lemma} The following properties hold
	\begin{itemize}
		\item[(\textit{i})] The relation $\prec$ is  ordering relation.
		\item[(\textit{ii})]  For every odd natural number $n$ and for every natural number $i>0$ the relation $n \prec S_{i}(n)$ holds.
		\item[(\textit{iii})]  For every  odd, natural number $n$ the relation $P(n)\prec n$ holds.
		\item[(\textit{iv})]  The function  $C(n)$   defines a permutation of the set $N$ of natural numbers.  
		\item[(\textit{v})] The graph $D=\bigcup\limits_{x=0}^{\infty}\,D_{c}$ is equal to the graph $\mathcal{G}$, i.e. $D= \mathcal{G}$. Therefore the graph $\mathcal{G}$ is connected.
	\end{itemize}
\end{lemma}
\begin{proof}
	A straightforward  proof of properties $(i)-(iv)$is left to the reader.  
	We shall prove the that the graph $\mathcal{G} $ is connected. \\
Suppose, on the contrary, that the graph G is not connected.
Hence there is a number $c$ such that the tree $D_{c+1}$ contains an element $n_{0}$ such that the element $n_{0}$ is not a successor of any element $n \in D_{c}$  (i.e.  for every number $n \in D_{c}$ and for every number $i$, the element $n_{0} \neq S_{i}(n)$).
On the other hand the predecessor $m=P(n_{0})$ of the number $n_{0}$ exists and $m \in D_{c}$	.
By the lemma \ref{zbr}(\ref{rwn14}) we have that there exists a number $j$ such that $S_{j}(m)=n_{0}$. This is inconsistent with our assumption. Therefore the graphs $D$ and $\mathcal{G}$ are equal,  $D=\mathcal{G}$.	
	\end{proof}
From this and the lemma \ref{acyl} we obtain
\begin{lemma}\label{Gtree}
	\color{white} blah \color{black}
\begin{enumerate}
	\item [(\textit{i})]	The graph $\mathcal{G}$ is a tree.
	\item [(\textit{ii})]  The graph $\mathcal{HC}$ is a tree.
\end{enumerate}
\end{lemma}

\section{Main lemma}\label{mle}

\begin{lemma}\label{MnLm}
	Let $r$ be a natural number. Every formula of the scheme \ref{MainLm} is a theorem of algorithmic theory $\mathcal{ATN}$.

\begin{small}
	\begin{equation}\label{MainLm}\tag{$ML_{r}$}
	\mathcal{ATN} \vdash
		(n=r) \implies
							 \{x:=0\} \bigcup\{x:=x+1\}\,\left( 
	\begin{array}{l} n\cdot 3^{x}+y=2^{z}  \land  \\ 
		\left| 	\begin{array}{l}
		\left( y=	\sum\limits_{j=0}^{x-1}\,3^{x-1-j}\cdot 2^{\sum_{p=0}^{j}k_{p}}  \right) \land \  
		\left(z= \sum\limits_{p=0}^{x}k_{p}\right) \land \\
		\ \underbrace{\left( \begin{array}{l}
			 (k_{0}=\kappa(n)\land m_{0}=\rho(n)  ) \land  \\ 
		\left( \bigwedge\limits_{l=0}^{x-1} \left( \begin{array}{l}
			 k_{l+1}=\kappa(3m_{l}+1) \land \\ m_{l+1}=\rho(3m_{l}+1)
		\end{array}
		\right)  \right)\end{array} \right) }_{\text{the length of computation }\ =x}   
	\end{array}\right| \end{array}\right) 	 	\quad
\end{equation}	
\end{small}
\end{lemma}
The lemma states that for every natural number $r$, the $3n+1$ computation terminates. More exacly, it states that the set $St_{0}$ of formulas of the form \ref{MainLm}  is recursive and it consists of theorems of $\mathcal{ATN}$ theory.
\begin{proof}
	The proof follows from the lemma \ref{Gtree} and the fact that the nodes of the tree $\mathcal{G}$ are standard, odd natural numbers. Hence no infinite computation exists.
\end{proof}

	\begin{footnotesize}
		\begin{example}\label{ex4.1M}
		$	\begin{array}{|r|l}
			n=5 \implies n\cdot 3^{\mathbf{1}} +1=16 \\
			n=67\implies n\cdot 3^{\mathbf{8}}+84701=2^{19} & \\
			n=5461\implies n\cdot 3^{1} +1=2^{14}		\\
		n=	227\implies n\cdot 3^{2}+5=2^{11} & \colorbox{green!10}{$\quad (\overbrace{(\overbrace{(227*3 +1)\ \div\ 2^{1})}^{P(227)=341}*3+1)\, \div\,2^{10}}^{P(P(227))}= 1$} 
		\end{array} $
	\end{example}
	\end{footnotesize}
	
	\begin{corollary}
		Let $St_{0}$ be the set of sentences of the form \eqref{s0}	
	\begin{enumerate}
		\item [(\textit{i})] every sentence of the set  $St_{0}$ is a theorem of $\mathcal{T}'$ theory (i.e. of an inessential extension of Presburger theory) and hence it is the theorem of $\mathcal{ATN}$ algorithmic theory of natural numbrs as well.
		\item [(\textit{ii})] the set $St_{0}$ is a recursive set.
	\end{enumerate}	
	\end{corollary}
	
	Let $r>0$ be a natural number, let $i_{r}$ be the number such that the formula \ref{MainLm} is valid in the structure $\mathfrak{N}$. 
	\begin{itemize}				
		\item 	The set of equations $St_{0}$ that contains all equalities of the  form \eqref{s0}  is a recursive set.
		\begin{equation} \label{s0}\tag{$S_{0}$}
			St_{0}\stackrel{df}{=} \left\lbrace \boxed{\ \ r\cdot 3^{i_{r}}+\underbrace{{\sum\limits_{j=0}^{i_{r}-1}3^{i_{r}-1-j}\cdot 2^{\sum_{p=0}^{j} k_{p}(r)}}}_{y} = 2^{ \underbrace{{\sum_{j=0}^{i_{r}} k_{j}(r) }}_{z} }}\right\rbrace_{r=1}^{\infty} 
		\end{equation}    
		\item Every element  of the set $St_{0}$  is a theorem of the theory  $\mathcal{ATP}$ (Preburer's arithmetic). 	
	\end{itemize} \color{black}
 	The lemma states that 1\textdegree) the infinite set $St_{0}$ of expressions is accompanied  by an algorithm that  decides whether a given formula $\varphi$ belongs to it, $\varphi \in St_{0}$ and 2\textdegree)  every element $\varphi$ of the set $St_{0}$ is a theorem of (inessentially etended) first-order theory of addition of natural numbers. Note, the proof of statement $\varphi$ is done by performing the computation of algorithm.
 
\begin{corollary}

		In other words, every number $r$ can be presented in the following form
		\begin{align}
			\frac{\frac{\frac{\frac{\frac{\frac{2^{k_{i_{r}}}-1 }{3}\cdot 2^{k_{i_{r}-1}}  -1}{3} \cdot 2^{k_{i_{r}-2}} -1}{3}    \cdot 2^{k_{i_{r}-3}} -1 }{3}  \cdots}{\cdots}\cdot 2^{k_{0}} -1}{3} = r   		
		\intertext{or in yet another  form }
			 r=S_{k_{0}}\Biggl(S_{k_{1}}\biggl(\cdots \Bigl(S_{k_{i_{r}-1}}\bigl( S_{k_{i_{r}}}(1)\bigr)\Bigr) \cdots\biggr)\Biggr) \\
	\end{align}
\end{corollary}
\begin{example}
	$\begin{array}{|c|}
		227=S_{1}(S_{5}(1))  \\
		61=S_{2}(S_{1}^2(S_{3}(S_{2}(1))))  \\
		27=S_{1}^{17}(S_{2}(S_{1}^{10}(S_{2}(S_{1}^{3}(S_{2}(S_{1}^{2}(S_{2}^{2}(S_{1}^{2}(S_{3}(S_{2}(1)) ))) ) ) ) )) )
	\end{array}$
\end{example}

\subsection{The limitations of the Main lemma}
\quad  \\
Does the lemma \ref{MnLm} solve the problem stated by Lothar Cllatz? \\
It seems so, since for every natural number $ r $ the computation of the algorithm $ Cl $ in the structure $ \mathfrak{N} $ is finite. \medskip\\
It turns out that the criterion of termination of computation  in its form \ref{twef} is \textit{falsifiable} in any non-standard model of Peaono's arithmetic in spite of being valid in the standard structure of natural numbers $\mathfrak{N}$.  
.. \\
Two remarks are in order.
\begin{remark}
	For every quadruple $\langle n,x,y,z \rangle$ of natural numbers the equation $n\cdot 3^{x}+y=2^{z}$ can be treated as an abbreviation of a formula of Presburger's arithmic, c.f. subsection \ref{prar} and example \ref{ex4.1M}.
\end{remark} \smallskip
 This  shows  that the sentence \ref{twef} is not a theorem. For the details consult the subsection  \ref{nst}.\\

\subsection{The positive consequnces of the main lemma. }\quad  \\
Let $\Gamma_{3}$ and $\Delta_{3}$ be parts of the program \ref{Gr3}.  \\
For every natural number $ r>0 $,  here exists the number $ i(r) $such that 
the implication   of the following form \ref{rTheta} 
\begin{equation}\label{rTheta}
n=r	\implies \left\lbrace \Gamma_{3}\,; (\mathbf{if}\ m \neq 1\ \mathbf{then}\ \Delta_{3}\ \mathbf{fi})^{i(r)} \right\rbrace (m=1) 
\end{equation} 
is a theorem of the algorithmic  theory $ \mathcal{ATN} $, c.f. the subsection \ref{ATN} on page \pageref{ATN}. 
\begin{equation}\label{rTheta2}
\mathcal{ATN} \vdash 	n=r\implies \left\lbrace \Gamma_{3}\,; (\mathbf{if}\ m \neq 1\ \mathbf{then}\ \Delta_{3}\ \mathbf{fi})^{i(r)} \right\rbrace (m=1) 
\end{equation} 
Note, that the formula \ref{rTheta} can be replaced by an equivalent formula (provided that the variable $l$ does not appear in  $\Delta_{3}$)
\begin{equation}\label{rTheta2}
	n=\underline{r}	\implies \left\lbrace \Gamma_{3}\,; \mathbf{for}\ l:= 1\ \mathbf{to}\ i(r)\   \mathbf{do}\ \Delta_{3}\ \mathbf{od} \right\rbrace (m=1) .
\end{equation}
\begin{example} 
\begin{small}
	$\begin{array}{|c|}
	\mathfrak{N} \models	n=67 \implies \left\lbrace m:=n;\,\text{for}\ l:=1\ \text{to}\ 8\  \text{do}\ m:=\frac{3*m+1}{2^{\kappa(3*m+1)}}    \ \text{od}\right\rbrace (m=1)  \\
	\mathfrak{N} \models n=373 \implies	\{m:=n;\,\text{for}\ l:=1\ \text{to}\ 4\  \text{do}\ m:=\rho(3*m+1)     \ \text{od}\}(m=1) \\
	\mathfrak{N} \models n=1367 \implies 	\{m:=n;\,\text{for}\ l:=1\ \text{to}\ 11\  \text{do}\ m:=\rho(3*m+1)     \ \text{od}\}(m=1) 
\end{array}$
\end{small}
\end{example}

\section{Proof of the Collatz conjecture}\label{dowd}	
Note, the formula \eqref{thM} expresses the conjecture of Collatz.
Look at the implication in the formula  \eqref{thM}. The \textit{antecedent} of this implication says: "\textit{$n$ is a natural number}" for it has the value $\mathbf{\mathbb{T}}$ (i.e. true) iff   $n=1 \lor n=2 \lor n=3 \lor \dots$. Similarly, the \textit{consequent} of the implication takes the value $\mathbf{\mathbb{T}}$ iff the program in the consequent terminates and the final value of the variable $m$ is 1.\\
\begin{theorem}[Main]\label{3nn+1}	
	The following formula {\tiny \eqref{thM}} is a theorem of formalized algorithmic theory $\mathcal{ATN}$
{\scriptsize 	\begin{equation}\label{thM} \tag{Main Thm}  
		\mathcal{ATN}\ \vdash \ \ 		
\underbrace{	\forall_{n \neq 0}	\ 
				\boxed{\left\lbrace \begin{array}{l}
						q:=1; \\
						\mathbf{while}\ n \neq q \ \mathbf{do}  \\
						\quad q:=q+1  \\
						\mathbf{od} 	 		
					\end{array}\right\rbrace(n=q)   }  }_{{\mathbb{FOR ALL}\ n\vert \ \ \ \ \ \ \ \mathbb{IF}\ n \ is\ a\ natural\ number}}
			\implies 
			\underbrace{\boxed{\left\lbrace \begin{array}{l}
						m:=n \div 2^{\kappa(n)}; \\
						\mathbf{while}\ m \neq 1 \ \mathbf{do}  \\
						\quad m:=3 \cdot m +1;  \\
						\quad m:= m \div 2^{\kappa(m)}  \\
						\mathbf{od} 	 		
					\end{array}\right\rbrace(m=1)}  }_{{\mathbb{THEN}\ the\ computation\   for\ n\ is\ finite\ \mathbb{FI}}}   					
	\end{equation}   } 
\end{theorem}
The theorem reads as follow, ({\footnotesize the program in the consequent of the implication is abbreviated as $Gr$}): \\
\textit{ for every $n$, if $ n \neq 0$ is a natural number, \\
	\hspace*{1.6cm} then the computation of the program  $ Gr $ is finite and final value of the variable  $ m=1 $.} \bigskip \\
The idea of the proof can be explained as follow:
\begin{enumerate}
	\item We know that the set $St_{0}$, see  page[\pageref{s0},\eqref{s0}], is recursive and consists of theorems of elementary theory of addition of natural numbers (i.e. the  Presburger's arithnetic with helpful but inessential extensions), and hence the  elements of the set $St_{0}$ are theorems  of the algorithmic theory of numbers $\mathcal{ATN}$. 
	\item We transform this set to three consecutive sets $St_{0} \rightarrow St_{1} \rightarrow St_{2} \rightarrow St_{3}$ in such a way that every set  $St_{i}, i=1,2,3$ is recursive and consists of theorems of algorithmic theory of natural numbers.  Every formula $\varphi$ of the set $St_{i+1}$ has a proof from some formula $\psi \in St_{i} $. 
	\item In the last step we apply the inference rule $R_{3}$, see page \pageref{Rules}, of infinitely many premises. \bigskip \\
\end{enumerate}
We use a couple of denotations.
Let $ r $ be a natural number.  \\ 
The sign $ \Gamma $ abbreviates the program   $\left\lbrace    m:= n \div 2^{\kappa(n)}; \right\rbrace $.  \\ 
The sign $ \Delta $ denotes the program $\left\lbrace  m:=3\cdot m +1;\ m:= m \div 2^{\kappa(m)} \right\rbrace $ \  or if you like 
 \[ \left\lbrace  m:=3\cdot m +1;\ \mathbf{while}\ even(n)\    \mathbf{do}\ n:=n \div 2 \ \mathbf{od} \right\rbrace  \]
\begin{proof}   
We are resuming at the main lemma \ref{MnLm}. For every natural number  $ r $ exists  number  $ i_{r} $ such that $ m_{i_{r}}=1 $.  \\

Consider the set $ St_{1} $ of formulas such that for every natural number $ r $ the formula of the scheme \eqref{S1z} 
	\begin{equation}\label{S1z}
	 	   n=\underline{r} \implies \left\lbrace \Gamma\right\rbrace  \left\lbrace \mathbf{if}\ m \neq 1 \ \mathbf{then}\ \Delta \ \mathbf{fi}\right\rbrace^{i(\underline{r})}  (m=1)      
\end{equation}
is contained in $ St_{1} $.  
\begin{lemma}
	The following  conditions hold
	\begin{itemize}
		\item[\textit{(i)}]  The set $St_{1}$ is a recursive set.
		\item[\textit{(ii)}] For every formula $\psi \in St_{1}$ there is a formula $\varphi \in St_{0}$ such rhat the equivalence $\varphi \equiv \psi$ is a theorem of program calculus. 
	\end{itemize}
\end{lemma}
	\begin{equation}\label{S11z}
	St_{1}\stackrel{df}{=}\left\lbrace    n=\underline{r} \implies \left\lbrace \Gamma\right\rbrace  \left\lbrace \mathbf{if}\ m \neq 1 \ \mathbf{then}\ \Delta \ \mathbf{fi}  \right\rbrace^{i(\underline{r})}  (m=1)     \right\rbrace   _{r=1}^{\infty}
\end{equation} 

			 In the proof of the lemma, we use the following axiom \eqref{asig} of assignment instruction 
			\begin{equation}
				\label{asig}  \tag{$Ax_{18}$} 	\fcolorbox{red}{red!10}{						 	
				$	\left\lbrace  x := \tau \right\rbrace \alpha(x)  \Leftrightarrow \alpha(x/\tau)	$}  
			\end{equation}
			\begin{flushright}
				{\scriptsize  On the lefthandside is an algorithmic formula where $\{x:=\tau\}$ is assignment instruction, and  $\alpha(x)$ is a formula.  On the righthandside is the formula that arises from $\alpha(x)$ by replacing all free occurencies of the variable  $x$ in $\alpha$ by expression $\tau$. }
			\end{flushright}
			and the axiom $Ax_{20}$ of conditional instruction \textbf{if}.  
				\begin{equation}
				\label{if1}  \tag{$Ax_{20}$} 	\fcolorbox{red}{red!10}{						 	
					$	\mathbf{if} \ \gamma \ \mathbf{then}\ K\ \mathbf{else}\ M\ \mathbf{fi}\,\alpha \Leftrightarrow 
					((\gamma \land K\alpha) \lor (\lnot M\alpha) ) 	$}  
			\end{equation}
Next, we consider the set  $ St_{2} $ that contain all  formulas of the scheme shown in the equation \ref{S2z} and only such  formulas.
	\begin{equation}\label{S2z}
		St_{2}\stackrel{df}{=}	\biggl\{  n=\underline{r} \implies \left\lbrace \Gamma;\  \mathbf{while}\ m \neq 1 \ \mathbf{do}\ \Delta \ \mathbf{od}\right\rbrace  (m=1)   \biggr\}_{r=1}^{\infty}
	\end{equation}
Our next observation says: 
     \begin{itemize}
     	\item[T2)]  Every formula  $ \phi $ from the set $ St_{2} $  is a theorem of  algorithmic theory  $ \mathcal{ATN} $ of natural numbers.  
     \end{itemize}
Each formula of the set $ St_{2}$ is proved from  a corresponding formula $ \varphi $ in the set $ St_{1} $ using  the following theorem \eqref{ThIf} of calculus of programs.  
	{\footnotesize \[\label{ThIf}  \tag{ThIf} 	\fcolorbox{red}{red!10}{$						 	
				\left\lbrace M;\mathbf{if}\ \gamma\ \mathbf{then}\ K\ \mathbf{fi}^{i}\right\rbrace (\alpha \land \lnot \gamma) \implies
				\left\lbrace      \begin{array}{l}
				 M;\ \  
				\mathbf{while}\ 
				\gamma\ \
				 \mathbf{do}\ K\ \mathbf{od} \end{array} \right\rbrace  (\alpha \land \lnot \gamma) $
	}\]  }
Hence the thesis   T2) is justified.  \\
Therefore, the set 
 $ St_{3} $ that consists of all formulas of the scheme  \eqref{S3z}.
	{\scriptsize 	\begin{equation}\label{S3z}
			St_{3}\stackrel{df}{=}	\left\{ \left\lbrace  \begin{array}{l}
				q:=0; \\
				\left\lbrace 
				\begin{array}{l}		    	
					q:=q+1 \\
				\end{array}
				\right\rbrace^{\underline{r}}
			\end{array} 
			\right\rbrace  
			(n=q)
			\implies 
			\left\lbrace \begin{array}{l}
				\Gamma;  \\
				\mathbf{while}\ m \neq 1  \\
				\mathbf{do}\\ \quad	 \Delta \\  \mathbf{od}
			\end{array}\right\rbrace (m=1) \right\}_{r=1}^{\infty}   
	\end{equation}     }
and no other formulas is the set of theorems of the theory $ \mathcal{ATN}  $	.  \\
Now, we apply the following inference rule
  $ R_{3} $ of calculus of programs
	\[\label{infR3}\tag{$ R_{3} $}  	\fcolorbox{red}{red!10}{$\dfrac{\bigl\{M; \mathbf{if}\ \gamma\ \mathbf{then}\ K\ \mathbf{fi}^{i}\bigr\}(\lnot \gamma \land\alpha) \implies \beta\}_{i=0}^{\infty}} 
		{\bigl\{M; \mathbf{while}\ \gamma\ \mathbf{do}\ K\ \mathbf{od}\bigr\}(\lnot \gamma \land\alpha)\implies \beta \}  } $} \]
to obtain the  theorem \ref{thM} of the  algorithmic theory   $ \mathcal{ATN} $ of natural numbers.
\end{proof}
\textsc{Comment}, the proof of the theorem \ref{thM} is an infinite tree, all of its branches are finite, all leafs are axioms (or some eaerlier proved theorems). Obviously,such  a proof can not be written in finite time. Instead, we have proved that the proof exists. For the definition of proof in the calculus of programs consult \cite{al:gm:as} Definition II.5.2, p.58.   \\
\begin{figure}
		\begin{normalsize}
			\begin{tikzcd}
				St_{0}\stackrel{df}{=}	\left\{\boxed{{n=\underline{r} \implies 
						\left( \overbrace{ \begin{array}{l}
								\bigl(r\cdot 3^{i_{r}}+y(i_{r})=2^{z_{i_{r}}}\bigr) \land \\
								\bigl(y(i_{r})=\sum\limits_{j=0}^{i_{r}-1}3^{i_{r}-1-j}\cdot 2^{\sum_{p=0}^{j} k_{p}(r)}\bigr) \land \\
								\bigl(z_{i_{r}}=\sum_{j=0}^{i_{r}}k_{j}(r) \bigr)
						\end{array} }^{\Psi(r,i_{r})}\right)  } 	}  
				\right\}_{r=1}^{\infty} \arrow[d,orange,very thick,"\mathrm{Lemma\ 5.1}"'] \\  
				St_{1}\stackrel{df}{=}	\left\{\boxed{{n=\underline{r} \implies \left\lbrace \Gamma\right\rbrace  \left\lbrace \mathbf{if}\ m \neq 1 \ \mathbf{then}\ \Delta \ \mathbf{fi}\right\rbrace^{i(\underline{r})}  (m=1)} } \right\}_{r=1}^{\infty} \arrow[d,brown,very thick]   
				\\
				St_{2}\stackrel{df}{=}	\biggl\{  \boxed{n=\underline{r} \implies \left\lbrace \Gamma;\  \mathbf{while}\ m \neq 1 \ \mathbf{do}\ \Delta \ \mathbf{od}\right\rbrace  (m=1) }  \biggr\}_{r=1}^{\infty}  \arrow[d,orange,very thick," r=0+\underbrace{1+1+\cdots+1}_{r\ \times}"']   
				\\\
				St_{3}\stackrel{df}{=}	\left\{ \boxed{\left\lbrace  \begin{array}{l}
						q:=0; \\
						\left\lbrace 
						\begin{array}{l}		    	
							q:=q+1 \\
						\end{array}
						\right\rbrace^{\underline{r}}
					\end{array} 
					\right\rbrace  
					(n=q)
					\implies 
					\left\lbrace \begin{array}{l}
						\Gamma;  \\
						\mathbf{while}\ m \neq 1  \arrow[d,orange,very thick,"\mathrm{\ \ \ \ applying\ the\ rule\ R3\ }"]  \\
						\mathbf{do}\	 \Delta \  \mathbf{od}
					\end{array}\right\rbrace (m=1) } \right\}_{r=1}^{\infty} \arrow[d,orange,very thick]   \\
				\forall_{n \neq 0}	\underbrace{\boxed{\left\lbrace \begin{array}{l}
							q:=1; \\
							\mathbf{while}\ n \neq q \ \mathbf{do}  \\
							\quad q:=q+1  \\
							\mathbf{od} 	 		
						\end{array}\right\rbrace(n=q)   }  }_{\color{black}{\mathbb{FOR ALL}\ n,\ \mathbb{IF}\ n \ is\ a\ natural\ number}}
				\implies 
				\underbrace{\boxed{\left\lbrace \begin{array}{l}
							m:=n \div 2^{\kappa(n)}; \\
							\mathbf{while}\ m \neq 1 \ \mathbf{do}  \\
							\quad m:=3 \cdot m +1;  \\
							\quad m:= m \div 2^{\kappa(m)}  \\
							\mathbf{od} 	 		
						\end{array}\right\rbrace(m=1)}  }_{\color{black}{\mathbb{THEN}\ the\ computation\   for\ n\ is\ finite\ \mathbb{FI}}}  \arrow[d,orange,very thick]
					\\	
					\mathcal{ATN} \vdash \	\forall_{n \neq 0}	
					{\boxed{\left\lbrace \begin{array}{l}
								m:=n \div 2^{\kappa(n)}; \\
								\mathbf{while}\ m \neq 1 \ \mathbf{do}  \\
								\quad m:=3 \cdot m +1;  \\
								\quad m:= m \div 2^{\kappa(m)}  \\
								\mathbf{od} 	 		
							\end{array}\right\rbrace(m=1)}  } & 
				\end{tikzcd}
			\end{normalsize}
			\caption{The structure of the proof  of Main theorem}
			\label{prdowd}
	\end{figure}
\section{Final remarks.}	\label{FinRem}

Our message does not limit itself to the proof of Collatz theorem.  We believe that the algorithmic logic is a proper tool in developing algorithmics, discrete mathematics and software engineering.\\
We are  presenting a solid argument that the algorithmic  language of program calculus is indispensable for expressing the semantic properties of programs. Halting property of program, correctness property, axiomatic specification of data structure of natural numbers,etc.,   can not be expressed by (sets) of first-order formulas.   \\
We show the potential of calculus of programs as a tool for:
\begin{itemize}
	\item specification of semantical properties of software and  
	\item verification of software against some specifications.
\end{itemize}
We hope the reader will forgive us for a moment of insistence  (is it a propaganda?). \\
Calculus of programs $\mathcal{AL} $ is a handy tool.  
For  there are some good reasons to use the calculus of programs 
\begin{large}
	\begin{itemize}
		\item[(\textit{i})]  The language of calculus $\mathcal{AL} $ contains algorithms (programs) and \textit{algorithmic formulas}  besides terms and first-order formulas.
		\item[(\textit{ii})] Any semantical property of an algorithm can be \textit{expressed} by an appropriate algorithmic formula. Be it termination, correctness or other properties.
		\item[(\textit{iii})] Algorithmic formulas enable to create complete, categorical \textit{specifications} of data structures  in  the form of algorithmic theories.
		\item[(\textit{iv})] Calculus of programs $ \mathcal{AL}$ offers the \textit{complete} set of tools for proving theorems  of algorithmic theories.  
\end{itemize}   \end{large} \bigskip 
\paragraph*{Tigthening the gap}\quad  \\
Another contribution that this paper offers is an original way of proving theorems. Our proof of the Collatz conjecture is made in two stages. First, we prove that there are some recursive sets of formulas consisting of theorems of the $\mathcal{ATN}$ theory. Second, we reason on the sets of theorems much like one reasons on formulas and apply the infinitary rules of inference to achieve the goal. \medskip\\
The G\"{o}del's theorem reveals a big gap between the set of valid properties of natural numbers and the set of theorems of Peano's arithmetic. We are hoping that the calculus of programs (\textit{aka} algorithmic logic) will permit to solve many questions e.g. the Goldbach conjecture.
\noindent There are many questions that remain \underline{open}.  
\paragraph*{Efficiency problems}\quad  \\
From the papers of M. Presburger \cite{Presb, Stans} and D. Cooper \cite{Cooper72} one can infer     an estimation of  the pessimistic cost of a $3n+1$ computaion is  $O(2^{2^{2^{n}}}) $
 \\
 From the proof of lemma \ref{l4.3} one can deduce that the cost of the Collatz algorithm is polynomial. \\
We have the following conjecture and a remark
\begin{conjecture}
	For every natural number $n$ its $3n+1$ computation consists of no more than $2n$ multiplications by 3 and no more than $3n$ divisions by 2.
\end{conjecture}
\begin{remark}
	One needs not to execute $3n+1$ computation. Just accept 1 as the \textit{proven} result.
\end{remark}

\subsection*{Acknowlegments} 
Andrzej Szałas has shown to us the lacunes in our earlier proofs. Antek Ciaputa helped in calculations and drawing of Hotel Collatz. Hans Langmaack, Wiktor Dańko, Paweł Gburzyński and Marek Warpechowski sent a couple of useful  comments. \bigskip\\

\section{Suplements}
In this section, for the reader's convenience, we have collected some definitions, facts and statements. some useful theorems.
\begin{enumerate}
	\item[\ref{nst}] The Jaśkowski's model of Presburger arithmetic for the use of programmers.
	\item[\ref{prar}]  A short presentation of Presburger theory of addition.
	\item[\ref{AL}] A short presentation of calculus of programs \textit{aka} algorithmic logic $\mathcal{AL}$.
	\item[\ref{ATN}.] A short presentation of the formalized algorithmic theory of natural numbers $\mathcal{ATN}$.
\end{enumerate}
\subsection{A structure  with counterexamples} \label{nst}.
\begin{flushright}
	\textit{{\scriptsize where Collatz computations may be of infinite length}}
\end{flushright}
Here we present some facts  that are were discovered around 1929 by two students of Alfred Tarski: Mojżesz Presburger \cite{Presb,Stans} and Stanisław Jaśkowski \cite{AT}. These facts are  less known to the IT community.\\
A programmer  may doubt the importance of those facts.
Yet, it is worth considering, non-standard data structures do exist, and this fact has ramifications.
Strange as they seem, still it is worthwhile to be aware of their existence.  \\

Now, we will expose the algebraic structure $ \mathfrak {J} $, which is a model of the theory  $Ar$, i.e. all axioms of   theory $Ar$ are true in the structure $ \mathfrak {J} $. First we will describe this structure as mathematicians do, then we will write a class (i.e. a program module) implementing this structure. \medskip \\
\subsubsection*{Mathematical description of Jaśkowski's  structure}
$ \mathfrak {J} $ is an algebraic structure
\[\tag{NonStandard} \label{Nst}\mathfrak{J}=\langle M;\, \mbox{\textbf{\underline{0}}},\mbox{\textbf{\underline{1}}},\oplus;\,= \rangle\]
such that $M \subset \mathbb{C}$ is a set of complex numbers  $k+\imath w $, i.e. of pairs  $\langle k,w \rangle$, where element  $k \in \mathbb{Z}$ is an integer,  and element $w \in \mathbb{Q}^+ $ is a rational, non-negative number $w \geq 0 $  and the following requirements are  satisfied: 
\begin {enumerate}
\item [(\textit {i})] for each element $  k+\imath w  $ if $ w = 0 $ then $ k \geq 0 $, 
\item [(\textit {ii})] \textbf{\underline{0}} $\stackrel{df}{=}$ \colorbox{green!30}{ $\langle 0+\imath 0 \rangle $}, 
\item [(\textit {iii})] \textbf{\underline{1}} $\stackrel{df}{=}$ \colorbox{green!30}{ $\langle 1.+\imath  0\rangle $}, 
\item [(\textit {iv})]  the operation $\oplus $ of addition is determined as usual
\[ (k+\imath w)  \oplus  (k '+\imath w')  \stackrel {df} {=} ( k + k ')+\imath (w + w') . \]
\item [(\textit {v})]   the predicate = denotes as usual identity relation.
\end {enumerate}
\begin{lemma}
	The algebraic structure $\mathfrak{J} $  is a model of  first-order arithmetic of addition of natural numbers $\mathcal{T} $. 
\end{lemma}
The reader may  check that every axiom of the $\mathcal{T} $ theory (see  definition\ref{thT}, p.\pageref{thT}), is a  sentence true in the structure $ \mathfrak {J} $, cf. next subsection \ref{prar}. \smallskip\\
The substructure $ \mathfrak {N} \subset \mathfrak {J} $  composed of only those elements for which $ w = 0 $ is also a model of the theory $\mathcal{T} $.   \\
It is easy to remark that elements of the form $\langle k,0 \rangle$ may be identified with natural numbers $k$, $k \in N$. Have a look at table \ref{tabe}\medskip \\
The elements of the structure  $ \mathfrak {N} $ are called \textit{reachable}, for they enjoy the following algorithmic property
\[\forall_{n \in N}\, \{y:=\textbf{0}; \mathbf{while} \ y \neq n\ \mathbf{do}\ y:=y+\textbf{1} \ \mathbf{od}\}(y=n)\]
The structure $ \mathfrak {J} $  is not a model of the $ \mathcal {ATN} $, algorithmic theory of natural numbers, cf . subsection \ref{ATN}.
Elements of the structure $ \langle k, w \rangle $. such as $ w \neq \textbf{0} $ are \textit {unreachable}. i.e. for each element $x_0=\langle k, w \rangle$ such that $w \neq0$ the following condition holds   \[\lnot \{y:=\textbf{0}; \mathbf{while} \ y \neq x_0\ \mathbf{do}\ y:=y+\textbf{1} \ \mathbf{od}\}(y=x_0)\] 
The substructure $ \mathfrak {N} \subset \mathfrak {J} $  composed of only those elements for which $ w = 0 $ is a model of the theory $ \mathcal {ATN} $  c.f. subsection \ref{ATN}.   The elements of the structure $ \mathfrak {N} $ are called \textit {reachable}. A  theorem of the foundations of mathematics states:
\begin {lemma}
The structures $ \mathfrak {N} $ and $ \mathfrak {J} $ are not isomorphic. 
\end {lemma}
For the proof see \cite{art:ag}, p. 256.
As we will see in a moment, this fact is also important for IT specialists. \bigskip \\
An attempt to visualize structure $\mathfrak{M} $  is presented in the form of table \ref{tabe}.
The universe of the structure $\mathfrak{J}$ decomposes onto two disjoint subsets (one green and one red). 
Every element of the form $\langle k,0 \rangle$  (in this case $k>0$)  represents  the natural number $k$.  Such elements are called \textit{reachable} ones. 
Note, 
\begin{definition}
	An element $n$ is a standard natural number (i.e. is \textit{reachable} ) iff the   program of adding ones to initial zero terminates
	\[n \in N \stackrel{df}{\Leftrightarrow}
	\{q:=\textbf{0};\, \textbf{while}\ q \neq n\ \textbf{do}\ q:=q+\textbf{1}\ \textbf{od} \}(q=n) \]
	or, equivalently
	\[n \in N  \stackrel{df}{\Leftrightarrow}
	\{ q:=\textbf{0} \} \bigcup \{ \mathbf{ if}\  n \neq  q\ \mathbf{then }\  q:=q+\textbf{1}\ \mathbf{fi }\}(q=n)\]
\end{definition} \medskip
 
\begin{tiny}
	\begin{table}[h]
		\caption{Model  $\mathfrak{J}$ of Presburger arithmetic    
			\begin{tiny}
				consists of complex numbers $a+\imath\,b $ where $b \in Q^+$ and $a \in Z$,  additional condition: $b=0 \Rightarrow a \geq 0$. Definition of order $n >m \stackrel{df}{\equiv} \exists_{ u\neq 0}\, m+u=n  $. Invention of  S. Jaśkowski (1929).  \end{tiny} }\label {tabe}	
				\resizebox{\textwidth}{!}{ $\begin{array}{c|c}
				\text{\textsc{{Standard}}\   elements} & \text{\textit{Unreachable}\  (\textsc{{infinite}}) elements} \\  \hline\hline
			
				&\begin{array}{c} 
					
					\begin{array}{cccccccccc}
						   &&&&\cdots \\ \hline
						{ -\infty  \cdots} & -11+ \imath 2 & -10+\imath 2 &\cdots  & 0+\imath 2 & 1+\imath 2     & & {\cdots \infty }  \\\hline
						 &&&&\cdots  \\ \hline
					{ -\infty  \cdots}  & -11+\imath \frac{53}{47} & -10+\imath \frac{53}{47} & \cdots  & 0+\imath \frac{53}{47} & 1+\imath \frac{53}{47}     & & {\cdots \infty }  \\ \hline 
						&&&&\cdots \\   \hline
						{ -\infty  \cdots}  & -11+\imath \frac{28}{49} & -10+\imath \frac{28}{49} & \cdots  & 0+\imath \frac{28}{49} & 1+\imath \frac{28}{49}     & & {\cdots \infty }  \\  \hline
						&&&&\cdots \\  \hline
						{ -\infty  \cdots}  & -11+\imath \frac{3}{47} & -10+\imath \frac{3}{47} & \cdots  & 0+\imath \frac{3}{47} & 1+\imath \frac{3}{47}     & & {\cdots \infty }  \\  \hline
						&&&&\cdots  
					\end{array} 
				\end{array} \\
				\begin{array}{ccccc}
					0 & 1 & 2 & \cdots   
				\end{array} \\ \hline\hline  
			\end{array}$     }
	\end{table}     
	\end{tiny}   

Note that the subset that consists of all non-reachable elements is well separated from the subset of reachable elements. Namely,  every reachable natural number is  less that any unreachable one. Moreover, there is no least element in the set of unreachable elements. I.e. the principle of minimum does not hold in the structure $ \mathfrak{M}$. \\
Moreover, for every element $n$ its computation contains either only standard, reachable numbers or is composed of only unreachable elements.
This remark will be of use in our proof. 
\begin{remark} \label{ReGr}
	For every element $n$ the whole Collatz computation is either in green or in reed quadrant of the table \ref{tabe}.
\end{remark}
Elements of the structure $\mathfrak{M}$  are ordered as usual
\[ \forall_{x,y} \ x<y \stackrel{df}{=} \exists_{z \neq \textbf{0}} \  x+z=y . \]
Therefore, each reachable element  is smaller than every unreachable element. \\
The order defined in this way is the lexical order. (Given two elements $p$ and $q$, the element lying higher is bigger, if both are of the same height then the element lying on the right is bigger.)\\
The order type is $\omega + (\omega^*+\omega)\cdot\eta$.
\begin{remark}
	The subset of unreachable elements (red ones on the table \ref{tabe}) does not obey the principle of minimum. 
\end{remark}
The mathematically oriented reader will easily accept our arguments. For the informaticians we offer the same thought in the  software's disguise.
\subsubsection*{Definition in a programming language}
Perhaps you have already noticed that the $ \mathfrak {M} $  is a computable structure. The following is a class 
that implements the structure $ \mathfrak {M} $. The implementation uses the integer type, we do not introduce rational   numbers explicitly. \medskip \\
\begin{scriptsize} \textsf{\begin{tabular}{l} \hline
			unit StrukturaM: class; \\
			\quad unit Elm: class(k,li,mia: integer); \\
			\quad begin \\
			\qquad if mia=0 or li * mia <0  or li=0 and k<0 then raise Error fi;  \\
			\quad end Elm; \smallskip\\
			\quad add: function(x,y:Elm): Elm; \\
			\quad \  begin \ 
			 result := new Elm(x.k+y.k,  x.li*y.mia+x.mia*y.li, x.mia*y.mia  )  \ 
		  end add; \smallskip\\
			\quad unit one : function:Elm;\ \ begin \ \ result:= new Elm(1,0,2)\ \ end one; \smallskip\\
			\quad unit zero : function:Elm;\ \ begin \ \ result:= new Elm(0,0,2)\ \ end zero; \smallskip\\
			\quad unit eq: function(x,y:Elm): Boolean;  \\
			 \qquad  begin \ \
			   result := (x.k=y.k) and (x.li*y.mia=x.mia*y.li ) \ \
			   end eq; \\
			end StrukturaM \ \\ \hline
\end{tabular}}  \end{scriptsize}
\medskip \\
The following lemma expresses the correctness of the implementation with respect to the axioms of Presburger arithmetic $ \mathcal{AP} $ (c.f. subsection \ref{prar}) treated as a specification of a class (i.e. a module of program).
\begin {lemma}
The structure $\mathfrak{E}=\langle E, add, zero, one, eq \rangle$ composed of the  set $E=\{o\  \mbox{object}: o\ in Elm \} $  of   objects of class Elm  with the \textit{add} operation is a model of the $ \mathcal{AP} $ theory, 
\[ \mathfrak{E} \models  \mathcal{AP}  \]
\end {lemma}
\subsubsection*{Infinite Collatz algorithm computation}
How to execute the Collatz algorithm in StructuraM? It's easy. \medskip \\
\begin{scriptsize} \textsf{\begin{tabular}{l} \hline
			pref StrukturaM block \\
			\quad var n: Elm; \\
			\quad unit odd: function(x:Elm): Boolean; ... result:=(x.k mod 2)=1  ... end odd; \\
			\quad unit div2: function(x:Elm): Elm; ... \\
			\quad unit 3xp1: function(n: Elm): Elm; \dots result:=add(n,add(n,add(n,one))); \dots end 3xp1; \\
			begin \\
			\quad n:= new Elm(8,1,2); \\
			\quad Cl: \begin{tabular}{|l|}  \hline
				while \ not eq(n,one) \ do \\
				\quad if \ odd(n)\ then  \\
				\quad\  \ \  n:=3xp1(n)  \  else\ n:= div2(n) \\
				\quad   fi \\
				od \\ \hline
			\end{tabular}  \quad (* \textit{a version of algorithm Cl that uses class Elm}  *) \\
			end block; \\ \hline			
\end{tabular}}   \end{scriptsize} \medskip\\
 The reader may experiment with our program. Note, it is easy to adapt to your favorite syntax (programming language). \\
Below we present the computation of Collatz algorithm for
$n=\langle 8,\frac{1}{2} \rangle$.
{\scriptsize \[\langle 8,\frac{1}{2} \rangle,\,\langle 4,\frac{1}{4} \rangle, \langle 2,\frac{1}{8} \rangle,\,\langle 1,\frac{1}{16} \rangle, \,\langle 4,\frac{3}{16} \rangle,\,\langle 2,\frac{3}{32} \rangle,\,\langle 1,\frac{3}{64} \rangle,\,\langle 4,\frac{9}{64} \rangle,\,\langle 2,\frac{9}{128} \rangle, \cdots   \]}
Note, the computation of algorithm $Gr$ for the same argument, looks simpler
{\scriptsize  \[\langle 8,\frac{1}{2} \rangle,\,\langle 4,\frac{1}{4} \rangle, \langle 2,\frac{1}{8} \rangle,\,\langle 1,\frac{1}{16} \rangle, \,\langle 1,\frac{3}{64} \rangle,\,\langle 1,\frac{9}{256} \rangle,\, \cdots   \]}
None of the elements of the above sequence is a standard natural number. Each of them is unreachable.
It is worth looking at an example of another calculation.
Will something change when we assign \textsf{n} a different object? e.g. \textsf {n: = new Elm (19,2,10)}? \\
\begin{scriptsize}\[ \begin{array}{l}\langle 19,\frac{10}{2} \rangle,\,\langle 58,\frac{30}{2} \rangle, \langle 29,\frac{30}{4} \rangle,\,\langle 88,\frac{90}{4} \rangle, \,\langle 44,\frac{90}{8} \rangle,\,\langle 22,\frac{90}{16} \rangle,\,\langle 11,\frac{90}{32} \rangle,\,\langle 34,\frac{270}{32} \rangle,\,\langle 17,\frac{270}{64} \rangle,  \\
		\langle 52,\frac{810}{64} \rangle,\langle 26,\frac{405}{64} \rangle,\langle 13,\frac{405}{128} \rangle,\langle 40,\frac{1215}{128} \rangle,\langle 20,\frac{1215}{256} \rangle,\langle 10,\frac{1215}{256} \rangle,\langle 5,\frac{1215}{512} \rangle,\langle 16,\frac{3645}{512} \rangle,\langle 8,\frac{3645}{1024} \rangle, \\
		\langle 4,\frac{3645}{2048} \rangle,\langle 2,\frac{3645}{4096} \rangle,\langle 1,\frac{3645}{8192} \rangle,\langle 4,\frac{3*3645}{8192} \rangle,\langle 2,\frac{3645*3}{2*8192} \rangle,\langle 1,\frac{3*3645}{4*8192} \rangle,\langle 4,\frac{9*3645}{4*8192} \rangle,\cdots   \end{array}\]  \end{scriptsize}
And one more computation.
\begin{scriptsize}  \[ \begin{array}{l}\langle 19,0 \rangle,\,\langle 58,0 \rangle, \langle 29,0 \rangle,\,\langle 88,0 \rangle, \,\langle 44,0 \rangle,\,\langle 22,0 \rangle,\,\langle 11,0 \rangle,\,\langle 34,0 \rangle,\,\langle 17,0 \rangle,  
		\langle 52,0 \rangle,\langle 26,0 \rangle, \\
		\langle 13,0 \rangle,\langle 40,0 \rangle,\langle 20,0 \rangle,\langle 10,0 \rangle,\langle 5,0 \rangle,\langle 16,0 \rangle,\langle 8,0 \rangle, \
		\langle 4,0 \rangle,\langle 2,0 \rangle,\langle 1,0 \rangle .\end{array}\]  \end{scriptsize}
\begin {corollary}
The structure $ \mathfrak {M} $, which we have described in two different ways, is the model of the $ \mathcal{AP}$ theory 

with the \underline{non-obvious} presence of unreachable elements in it.
\end {corollary}
\begin {corollary}
The halting property of the Collatz algorithm cannot be proved from the axioms of the $ \mathcal{T} $ theory, nor from the axioms of $ \mathcal{AP}$ theory.
\end {corollary}
\subsection{Presburger's arithmetic} \label{prar}\quad  \\
Presburger's arithmetic is another name of  elementary theory of natural numbers with addition. \\
We shall consider the following theory , cf. \cite{Presb},\cite{art:ag} p. 239 and following ones.
\begin{definition}\label{thT}
	Theory $\mathcal{T}=\langle \mathcal{L},\mathcal{C}, Ax \rangle$ is the system of three elements:\begin{itemize} 
		\item[$\mathcal{L}$] is a language of first-order.   The alphabet of this language consist of:  the set  $V$ of variables, symbols of operations: $0,S, + $, symbol of equality relation  $=$, symbols of logical functors and quantifiers,  auxiliary symbols as brackets. \\
		The set of well formed expressions is the union of te set   $T$ of terms and the set of formulas  $F$. \\
		The set  $T$ is the least set of expressions that contains the set  $V$  and constants 0 and 1 and closed with respect to the rules:  \begin{enumerate}
			\item[a)] if two expressions  $\tau_1$ and $\tau_2$ are terms, then  the expression $(\tau_1+\tau_2)$ is a term too.
			\item[b)]  if an expression $\tau$ is a term, then the expression $S(\tau)$ is a term too.
		\end{enumerate} 
		The set $F$ of formulas  is the least set of expressions that contains the equalities  (i.e. the expressions of the form $(\tau_1 =\tau_2)$) and closed with respect to the following formation rules: if expressions  $\alpha$ and  $\beta$ are formulas, then the aexpression of the form 
		\[(\alpha \lor \beta) ,\  (\alpha \land \beta),\  (\alpha \implies \beta) , \ \lnot\alpha\]
		are also formulas, moreover, the expressions of the form  
		\[\forall_x\,\alpha, \ \exists_x\,\alpha\]
		where  $x$ is a variable and  $\alpha$ is a formula, are formulas too. \\
		\item[$\mathcal{C}$] is  the operation of   consquence  determined by axioms of first-order logic and the inference rules of the logic, 
		\item[$Ax$] is the set of formulas listed below.  
	\end{itemize}
	\label{AX} 
	\begin{flalign} 
		\tag{a}  & \forall_x\ x+1 \neq 0 & \\
		\tag{b}  & \forall_x\, \forall_y\ x+1=y+1 \implies  x=y & \\
		\tag{c}  & \forall_{x}\ x+0=x & \\
		\tag{d}  & \forall_{x,y}\ (y+1)+x=(y+x)+1 & \\
		\tag{I}  &  \Phi(0)\land \forall_x\,[\Phi(x) \implies \Phi(x+1)]\implies \forall_x\Phi(x)   
	\end{flalign}     
	The expression $\Phi(x)$ may be replaced by any formula.   The result is an axiom of theory 
	This is the induction scheme.  \\
	We augment the set of axioms adding four axioms that define a coiple of useful notions.\\
	\begin{flalign*}
		\tag{e} & even(x) \stackrel{df}{\equiv} \exists_y\, x=y+y & \\
		\tag{o} & odd(x) \stackrel{df}{\equiv} \exists_y\, x=y+y+1 & \\
		\tag{D2} & x\, div\, 2 = y \equiv (x=y+y\, \lor\, x=y+y+1)  &\\
		\tag{3x} & 3x\stackrel{df}{=} x+x+x
	\end{flalign*} 
\end{definition}
The theory $\mathcal{T}'$ obtained in this way is a conservative extension of theory $\mathcal{T}$.\bigskip \\
Below we present another theory $\mathcal{AP}$ c.f. \cite{Presb}, we shall use two facts: 1) theory $\mathcal{AP}$ is complete and hence is decidable, 2) both theories are elementarily equivalent.
\begin{definition}
	Theory  $\mathcal{AP}=\langle \mathcal{L},\mathcal{C}, AxP \rangle$ is a system of three elements :\begin{itemize} 
		\item[$\mathcal{L}$] is a language of first-order. The alphabet of this language contains the set  $V$ of variables, symbols of functors : $0, + $, symbol of equality predicate $=$. \\ The set of well formed-expressions is the union of set of terms $T$  and set of formulas  $F$.
		The set of terms  $T$ is the least set of expressions that contains the set of variables  $V$  and the expression $0$ and closed with respect to the following two rules:   1) if two expressions  $\tau_1$ and $\tau_2$ are terms, then the expression  $(\tau_1+\tau_2)$ is also a term,  2) if the expression  $\tau$ is a term, then the expression  $S(\tau)$ is also a term.
		\item[$\mathcal{C}$] is the consequence operation determined by the axioms of predicate calculus and inference rules of first-order logic 
		\item[$AxP$] The set of axioms of the $\mathcal{AP}$ theory is listed below.
	\end{itemize}
	\label{AXMS} \begin{flalign*}  
		\tag{A}  & \forall_x\ x+1 \neq 0 & \\
		\tag{B}  & \forall_x\ x\neq 0 \implies \exists_y x=y+1 & \\
		\tag{C}  & \forall_{x,y}\ x+y=y+x & \\
		\tag{D}  & \forall_{x,y,z}\ x+(y+z)=(x+y)+z & \\
		\tag{E}  & \forall_{x,y,z}\ x+z=y+z \implies x=y & \\
		\tag{F}  & \forall_x\ x+0=x & \\
		\tag{G}  & \forall_{x,z}\ \exists_y\ (x=y+z \lor z=y+x) & \\
		\tag{H2} &  \forall_x\, \exists_y\ (x=y+y \lor x=y+y+1) & \\
		\tag{H3} &  \forall_x\, \exists_y\ (x=y+y+y \lor x=y+y+y+1 \lor x=y+y+y+1+1) & \\
		\intertext{. . . . . . . . . }
		\tag{Hk} &  \forall_x\, \exists_y\ \left( \begin{array}{l} x=\underbrace{y+y+\dots+y}_k \,\lor  \\ \ \ \ x=\underbrace{y+y+\dots+y}_k+1 \,\lor \\
			\ \ \ \ \ \  x=\underbrace{y+y+\dots+y}_k+\underbrace{1+1}_2 \,\lor  \\
			\qquad \dots \\
			\ \ \ \ \ \ \ \ \      x=\underbrace{y+y+\dots+y}_k+\underbrace{1+1+\dots+1}_{k-2}\, \lor \\
			\ \ \ \ \ \ \ \ \ \ \ \  x=\underbrace{y+y+\dots+y}_k+\underbrace{1+1+\dots+1}_{k-1} 
		\end{array} \right)   &  \\
		\intertext{\qquad \dots } \\  
	\end{flalign*}
\end{definition}
Let us recall a couple of useful  theorems   \smallskip\\
\textbf{F1}. Theory $\mathcal{T}$ is elementarily equivalent to the theory  $\mathcal{AP}$.\cite{Presb} \cite{Stans} \smallskip\\
\textbf{F2}. Theory $\mathcal{AP}$ is decidable. \cite{Presb}.   \smallskip\\
\textbf{F3}. The computational complexity of theory  $\mathcal{AP}$,   is triple exponential  $O(2^{2^{2^{cn}}})$. This result belongs to Oppen  \cite{oppen} and it can be traced in the Presburger's paper \cite{Presb}.  \smallskip\\
\textbf{F4}.  Theories  $\mathcal{T}$
and  $\mathcal{AP}$ have non-standard model, see section \ref{nst}, p. \pageref{nst}.  \smallskip \\

The reader may have doubts about the correctness of the formulas in our text that are using exponentiation.
Of course, exponentiation does not occur in Presburger arithmetic.  \\
We are offering the following explanations andcomments.
\begin{itemize}
	\item Let $n$ be a variable and $x$ be a natural number (e.g., 3). Then, the notation n . $n\cdot 3^4$ should be understood as an abbreviation for the expression $(n+n+ ....+n)$ in which the variable $n$ has been added 81 times.  
	\item Let $a,x,y,z$ be numerals or terms withou variables.  The formula $a\cdot 3^{x}+y=2^{z}$ comes as a neat abbreviation of some lengthy and obscure expression.
	\item  In the section \ref{dowd} we are operating within algorithmic theory $\mathcal{ATN}$. This theory is alike the Presburgers theory where the language admits algorithms (programs) as well formed expressions and instead of scheme of induction we have the axiom of reachability (R), \ on page \pageref{ATN}.  \\
	Consider the following program $P3(n,x)$ where $n$ and $x$ are some variables.
	\begin{equation*}P3(n,x) :	\left\lbrace \begin{array}{l} 
			q \leftarrow 0;\, w \leftarrow n;  \\
			\mathbf{while}\ q \neq x\   \mathbf{do} \\
			\quad w\leftarrow w+w+w;  \ \
			\quad q \leftarrow q+1  \\
			\mathbf{od} 	
		\end{array}  \right\rbrace 
	\end{equation*}
	If the initial value of the variable $n$ is a number $r$ and the initaila value of the variable $x$ is a number $a$, i.e. $val: \dfrac{n\ \vert \ x}{r\ \vert \ a}$, then the execution of the program $P3$ is finite and the final value of the variable $w= r^{a}$. 
\end{itemize}

	\subsection{An introduction to the calculus of programs $\mathcal{AL}$} \label{AL}\quad  \\
	For the convenience of the reader we cite the axioms and inference rules of calculus of programs i.e. algorithmic logic $\mathcal{AL}$.\\
	\textbf{Note}. \textit{Every axiom of algorihmic logic is a  tautology. \\
		Every inference rule of $\mathcal{AL}$ is sound.} \cite{al:gm:as}\bigskip\\
	\begin{normalsize}\textbf{Axioms}      \end{normalsize}  
		\begin{footnotesize}
 		\begin{alignat*}{4}
			\intertext{\textit{axioms of propositional calculus}}
			Ax_1) & \left(   \begin{array}{l}
				((\alpha\Rightarrow\beta
				)\Rightarrow((\beta\Rightarrow\delta)\Rightarrow 
				(\alpha
				\Rightarrow\delta))) 
			\end{array}\right)    & Ax_2)  &\ \  (\alpha\Rightarrow(\alpha
			\vee\beta)) \\
			Ax_3) & (\beta\Rightarrow(\alpha\vee\beta))  &       Ax_4) & ((\alpha\Rightarrow\delta)~\Rightarrow((\beta
			\Rightarrow\delta)~\Rightarrow  \\
			& & & ~((\alpha\vee\beta
			)\Rightarrow\delta)))  \\
			Ax_5)  & ((\alpha\wedge\beta	)\Rightarrow\alpha) & Ax_6) & ((\alpha\wedge\beta	)\Rightarrow\beta)   \\
			Ax_7) & \left(   \begin{array}{l}
				((\delta\Rightarrow\alpha
				)\Rightarrow((\delta\Rightarrow\beta)\Rightarrow  
				(\delta
				\Rightarrow(\alpha\wedge\beta)))) 
			\end{array}\right)  & Ax_8) & ((\alpha
			\Rightarrow(\beta\Rightarrow\delta))\Leftrightarrow((\alpha\wedge\beta
			)\Rightarrow\delta))   \\	
			Ax_9) & ((\alpha\wedge\lnot\alpha) \Rightarrow\beta) & Ax_{10}) & ((\alpha\Rightarrow
			(\alpha\wedge\lnot\alpha))\Rightarrow\lnot\alpha)   \\
			Ax_{11}) & (\alpha\vee\lnot\alpha)      \\
			\intertext{\textit{axioms of predicate calculus}}
			Ax_{12})\  & ((\forall x)\alpha(x)\Rightarrow \alpha(x/\tau))) & Ax_{13})\  & (\forall x)\alpha(x)\Leftrightarrow\lnot
			(\exists x)\lnot\alpha(x)   \\
			\intertext{\textit{axioms of calculus of programs} }
			Ax_{14})\  & K((\exists x)\alpha
			(x))\Leftrightarrow(\exists y)(K\alpha(x/y)) & Ax_{15})\  & K(\alpha\vee\beta)\Leftrightarrow((K\alpha)\vee(K\beta))   \\
			Ax_{16})\  & K(\alpha\wedge\beta)\Leftrightarrow((K\alpha)\wedge(K\beta)) & Ax_{17})\  & K(\lnot\alpha)\Rightarrow\lnot(K\alpha)   \\
			Ax_{18})\  &\left(   \begin{array}{l}
				 ((x:=\tau)\gamma\Leftrightarrow 
				 (\gamma(x/\tau)\wedge
				 (x:=\tau)true)) 
			\end{array}\right)     
			& Ax_{19})\  & \mathbf{begin}\ K;M\ \mathbf{end} \alpha
			\Leftrightarrow K(M\ \alpha)   \\
			 Ax_{20})\ & \left(\begin{array}{l}
				 \mathbf{if} \ \gamma \ \mathbf{then}\ K\ \mathbf{else}\ M\ \mathbf{fi}\,\alpha \Leftrightarrow 
				 ((\gamma \land K\alpha) \lor (\lnot M\alpha) )
			\end{array} \right)  &&  \\
		  Ax_{21})\ &   \left(   \begin{array}{l}
			 \mathbf{while}\ 			 \gamma \ \mathbf{do}\  K\ \mathbf{od}\  \alpha\Leftrightarrow   \\
			 (\lnot\gamma\wedge\alpha)\lor  
			 \bigl(\gamma\wedge K \ \mathbf{while}\  \gamma \ \mathbf{do}\ K \ \mathbf{od} (\lnot\gamma\wedge\alpha))\bigr) 
			\end{array}\right)    && \\
			Ax_{22})\  & {\bigcap K\alpha\Leftrightarrow				(\alpha\wedge(K\bigcap K\alpha))} & Ax_{23})\  & {\bigcup K\alpha				\equiv(\alpha\vee(K\bigcup K\alpha))}     		
		\end{alignat*}  
		\end{footnotesize}
	\begin{normalsize}\textbf{Inference rules}\end{normalsize} \label{Rules}
\begin{alignat*}{4} 
	\intertext{\textit{propositional\  calculus}}
	R_1\   & \dfrac{\alpha ,(\alpha \Rightarrow \beta )}{\beta }  & \hspace*{3cm}\\
	\intertext{\textit{predicate calculus}}
	R_{6}\  &  \dfrac{(\alpha (x)~\Rightarrow ~\beta )}{((\exists
		x)\alpha (x)~\Rightarrow ~\beta )}   &     R_{7}\  & \dfrac{(\beta ~\Rightarrow ~\alpha (x)}{(\beta
		\Rightarrow (\forall x)\alpha (x))}   \\
		\intertext{\textit{calculus\ of\ programs AL}}
		R_{2}\  & \dfrac{(\alpha \Rightarrow \beta )}{(K\alpha \Rightarrow K\beta ) }   &     R_{3}\  & \dfrac{\{s(\{\mathbf{if}\
			\gamma \ \mathbf{then}\ K\ \mathbf{fi}\}^{i} \ \alpha )\Rightarrow \beta \}_{i\in N}}{s(\{\mathbf{while}\
			\gamma \ \mathbf{do}\ K\ \mathbf{od}\} \ \alpha )\Rightarrow \beta }  \\
			\ \\
			R_{4}\  & \dfrac{\{(K^{i}\alpha \Rightarrow \beta )\}_{i\in N}}{(\bigcup				K\alpha \Rightarrow \beta )}   &        R_{5}\  & \dfrac{{(\alpha \Rightarrow K^{i}\beta )}_{i\in N}}{(\alpha				\Rightarrow \bigcap K\beta )} 
\end{alignat*}

	In the rules $R_6$ and $R_7$, it is assumed that $x$ is a variable which is not free in $\beta $, i.e. $x \notin  FV(\beta )$. The rules are known as the rule for
	introducing an existential quantifier into the antecedent of an implication
	and the rule for introducing a universal quantifier into the suc\-ces\-sor
	of an implication. The rules $R_4$ and $R_5$ are algorithmic counterparts of rules
	$R_6$ and $R_7$. They are of a different character, however, since their sets of
	premises are in\-finite. The rule $R_3$ for introducing a \textbf{while} into the
	antecedent of an implication of a similar nature. These three rules are
	called $\omega $-rules.
	The rule $R_{1}$ is known as \textit{modus ponens}, or the \textit{cut}-rule.
	In all the above schemes of axioms and inference rules, $\alpha $, $\beta $, 
	$\delta $ are arbi\-trary for\-mulas, $\gamma $ and $\gamma ^{\prime }$ are
	arbitrary open formulas, $\tau $ is an arbitrary term, $s$ is a finite
	se\-quence of assignment instructions, and $K$ and $M$ are arbitrary
	programs. \smallskip \\
	\begin{theorem}[\textit{theorem on  completeness of the calculus }   $\mathcal{AL}$].
		Let $\mathcal{T}=\langle \mathcal{L}, \mathcal{C}, \mathcal{A}x\rangle$ be a consistent algorithmic theory, let $\alpha \in \mathcal{L}$ be a formula. The following conditions are equivalent \begin{enumerate}
			\item[(\textit{i})] Formula $\alpha$ is a theorem of the theory T, \ $\alpha \in \mathcal{C}(\mathcal{A}x) $,
			\item[(\textit{ii})]  Formula $\alpha$ is valid in every model of the theory T, \ $ \mathcal{A}x \models \alpha$.
		\end{enumerate}
	\end{theorem}
	The proof may be found in \cite{al:gm:as} Theorem III.2.5 p.94.  

	\subsection{An outline of the algorithmic theory of natural  numbers $\mathcal{ATN}$} \label{ATN} \quad  \\
	The language of algorithmic theory of natural numbers $\mathcal{ATN}$ is very simple. Its alphabet contains one constant  0  \textit{zero }, one  one-argument functor $s$ and predicate = of equality. \\
	Axioms of  $\mathcal{ATN}$  were presented in the book \cite{al:gm:as}\\
	\begin{center}
		\begin{align*}
			(R)\  & \forall x \{ q:=0;\,\textbf{while}\ q \neq x\ \textbf{do}\ q:=s(q)\ \textbf{od}\}( q=x) \\ 
			(N)\  & \forall x \ s(x) \neq 0  \\ 
			 (J)\  & \forall x \forall y\, s(x)=s(y) \implies x=y   \\
				(D)\  & \forall x\,\forall y \left\{ \begin{array}{l} q:=0;w:=x;\\\mathbf{while}\ q \neq y\ \mathbf{do} \quad  q:=s(q)\ ; w:=s(w) \mathbf{od}  \end{array}  \right\} ( x+y=w)  
		\end{align*} 
	\end{center}
	
	The termination property of the program in $ (D) $is a theorem of  $\mathcal{ATN}$ theory as well as the formulas $x+0=x$ and $x+s(y)=s(x+y)$.  \\
 
	{\footnotesize Note, the following formulas are not theorems of first-order arithmetic of natural numbers (Peano's theory).}
	\begin{scriptsize} \begin{align}
			\label{l13}\mathcal{ATN}  & \vdash     {\exists_x\,\alpha(x) \Leftrightarrow   \{x:=0\}\bigcup\{x:=x+1\}\alpha(x)} \\
			\label{l14}\mathcal{ATN}  & \vdash     {\forall_x\,\alpha(x) \Leftrightarrow   \{x:=0\}\bigcap\{x:=x+1\}\alpha(x)} \\
			\intertext{The axiom of of Archimedes -- asserts the Archimedean property of natural numbers}
			\mathcal{ATN}  & \vdash0<x<y \implies \{a :=x;\, \mathbf{while} \ a<y \ \mathbf{do}\ a:=a+x\ \mathbf{od}\}(a\geq y)
			\intertext{Scheme of induction -- let $\alpha(x)$ be an algorithmic formula}  
			\label{l15}\mathcal{ATN}   & \vdash   \color{black}{\Bigl(\alpha(x/0) \land \forall_x\bigl(\alpha(x) \implies \alpha(x/s(x))\bigr)\Bigr)\implies \forall_x \alpha(x)}  \\
			\intertext{Correctness of Euclid's algorithm}   
			\label{l16}\mathcal{ATN}   & \vdash \left(\begin{array}{l}n_0>0 \land \\ m_0>0\end{array}\right)\implies\left\{\begin{array}{l}n:=n_0;\ m:=m_0; \\
				\mathbf{while}\,n \neq  m\, \mathbf{ do }\\ 
				\quad \mathbf{if }\,n>m \ \mathbf{then }\  \ n:=n\stackrel{.}{\_}m \\ 
				\qquad \mathbf{else }\ \ m:=m\stackrel{.}{\_}n \\ \quad \mathbf{fi}\\ \mathbf{ od} \end{array}\right\} (n=gcd(n_0,m_0)) \  
	\end{align}  \end{scriptsize}
	The theory $\mathcal{ATN}$  enjoys an important property of categoricity.
	\begin{theorem}[ \textit{meta}-theorem on categoricity of  $\mathcal{ATN}$ theory] Every model $\mathfrak {A}$ of the algorithmic theory of natural numbers is isomorphic to  the structure $\mathfrak{N}$. 
	\end{theorem}
	Compare this theorem with the fact that first-order theory of natural numbers has non-standard model, c.f. subsection \ref{nst}.   \bigskip \\
	As the simple corollary from the categoricity theorem we note \\
	\textit{every formula $\alpha$ valid in the standard model of natural numbers $\mathfrak{N}$ has a proof in the theory $\mathcal{ATN}$.} \\
	We can formulate this remark as follows: \textit{for every algorithmic formula $\alpha$, if the formula is valid in the structure $\mathfrak{N}$ then there is a finite proof of it or a finite proof that the finite prof exists or a finite proof of the existence of a proof thea t finite proof exists or \dots } \\
	
\end{document}